\documentclass[12pt]{article}

\usepackage{amsmath}
\usepackage{amssymb}
\usepackage{amsthm}
\usepackage{color}
\usepackage{graphicx}
\usepackage{color}

\def\be{\begin{eqnarray}}
\def\ee{\end{eqnarray}}

\def\b*{\begin{eqnarray*}}
\def\e*{\end{eqnarray*}}

\newtheorem{Theorem}{Theorem}[section]
\newtheorem{Definition}[Theorem]{Definition}
\newtheorem{Proposition}[Theorem]{Proposition}

\newtheorem{Assumption}[Theorem]{Assumption}

\newtheorem{Lemma}[Theorem]{Lemma}
\newtheorem{Corollary}[Theorem]{Corollary}
\newtheorem{Remark}[Theorem]{Remark}

\makeatletter \@addtoreset{equation}{section}

\newcommand*{\TitleFont}{
      \fontsize{19}{25}%
      \selectfont}

\usepackage{color}

\def\bit{\begin{itemize}}
\def\eit{\end{itemize}}

\def \E{\mathbb{E}}
\def \F{\mathbb{F}}
\def \P{\mathbb{P}}
\def \Ph{\widehat \mathbb{P}}
\def \Q{\mathbb{Q}}
\def \R{\mathbb{R}}

\def\Ac{{\cal A}}
\def\Bc{{\cal B}}

\def\Dc{{\cal D}}
\def\Ec{{\cal E}}
\def\Fc{{\cal F}}

\def\Ic{{\cal I}}
\def\Jc{{\cal J}}
\def\Kc{{\cal K}}
\def\Lc{{\cal L}}
\def\Mc{{\cal M}}
\def\Nc{{\cal N}}
\def\Oc{{\cal O}}
\def\Pc{{\cal P}}

\def\Pch{\widehat{\Pc}}

\def\Tc{{\cal T}}
\def\Uc{{\cal U}}

\def\Xc{{\cal X}}
\def\Yc{{\cal Y}}
\def\Zc{{\cal Z}}

\def\Rb{\overline{\R}}

\def\xb{\mathbf{x}}

\def\Fch{\widehat{\Fc}}

\def \Ph{\widehat{\P}}
\def \tauh{\hat{\tau}}

\def \I{{\bf I}}

\def \O{\Omega}
\def \Om{\Omega}

\def \Omh{\widehat{\Om}}
\def \Yh{\widehat{Y}}

\def \o{\omega}
\def \om{\o}
\def \omh{ \hat{\om}}
\def \wh{ \hat{\w}}

\def \eps{\varepsilon}

\def \Qh {\widehat{\Q}}
\def \w{\mathsf{w}}
\newcommand{\rmi}{{\rm (i)$\>\>$}}
\newcommand{\rmii}{{\rm (ii)$\>\>$}}
\newcommand{\rmiii}{{\rm (iii)$\>\,    \,$}}
\newcommand{\rmiv}{{\rm (iv)$\>\>\,$}}
\newcommand{\rmv}{{\rm (v)$\>\>\,$}}
\newcommand{\rmvi}{{\rm (vi)$\>\>\,$}}

\newcommand{\rma}{{\rm a)$\>\>$}}
\newcommand{\rmb}{{\rm b)$\>\>$}}
\newcommand{\rmc}{{\rm c)$\>\>$}}

\def\x{\times}
\def\ox{\otimes}
\def\1{{\bf 1}}
\def \proof{{\noindent \sc Proof. }}
\def\edoc{\end{document}}

\title{ \TitleFont Capacities, Measurable Selection \& Dynamic Programming\\
		Part I: Abstract Framework \footnote{Research supported by the chaire ``Risque Financier'' of the ``Fondation du Risque''.}}

\author{ El Karoui Nicole
\footnote{LPMA, UMR CNRS  6632, UPMC(ParisVI), and CMAP, UMR CNRS 7641, \'Ecole Polytechnique.}
\\
\and Tan Xiaolu
\footnote{CEREMADE, UMR CNRS 7534,  Universit\'e Paris Dauphine. }}

\date{\today}

\begin{document}

\maketitle

\abstract{ We give a brief presentation of the capacity theory and show how it derives naturally a measurable selection theorem following the approach of Dellacherie \cite{Dellacherie_1972}.
	Then we present the classical method to prove the dynamic programming of discrete time stochastic control problem, using measurable selection arguments.
	At last, we propose a continuous time extension, that is an abstract framework for the continuous time dynamic programming principle (DPP).
}
\vspace{1mm}


\vspace{1mm}


\section{Introduction}

	The capacity was first introduced by Choquet \cite{Choquet_1955} to derive an approximation property of Borel sets by its compact subsets in the real number space.
    It is then extended by himself in \cite{Choquet_1959} to an abstract form which generalizes the measures on a measurable space.
    Since then, the theory of capacities was used and developed in several ways, see e.g. Dellacherie \cite{Dellacherie_1972, Dellacherie_1972b}, etc.
    In particular, it gives a simple and straightforward proof of the measurable selection theorem, presented in Dellacherie \cite{Dellacherie_1972}, Dellacherie and Meyer \cite{Dellacherie_Meyer_1978}.

	As an important topic in the set theory, measurable selection theorem is applied in many fields, such as optimization problem, game theory, stochastic process, dynamic programming etc.
	Let $X$ and $Y$ be two spaces, $2^Y$ denote the collection of all subsets of $Y$, and $F:X \to 2^Y$ be a set-valued mapping.
	Then the question is whether there is a ``measurable'' mapping $f: X \to Y$ such that $f(x) \in F(x)$ for every $x \in X$.
	The selection theorems of Dubins and Savage \cite{Dubins_Savage}, Kurotowsky and Ryll-Nardzewski \cite{Kuro_Ryll_Nard} are usually given in this form.
	Another way to present the measurable selection theorem is to consider the product space $X \x Y$ with a subset $A \subseteq X \x Y$,
	and then to search for a ``measurable'' mapping $f:X \to Y$ such that its graph $[[f]] := \{ (x, f(x)) ~: x \in X \} \subseteq A$. Jankov-von Neumann's \cite{VonNeumann} theorem is given in this way.
	In this case, the graph set $[[f]]$ can be viewed as a section of the set $A$ in the product space $X\x Y$, and therefore it is also called the measurable section theorem.
	We refer to Parthasarathy \cite{Parthasarathy} and Srivastava \cite{Srivastava} for a detailed presentation of different selection theorems.

	The dynamic programming principle (DPP) is a principle which splits a global time optimization problem into a series of local time optimization problems in a recursive manner.
	It plays an essential role in studying the control problems, such as, to derive a computation algorithm, to obtain a viscosity solution characterization of the value function, etc.
	To derive a DPP of a stochastic control problem, it is also classical to use the measurable selection theorem as shown in Dellacherie \cite{Dellacherie_MaisonJeux} as well as in Bertsekas and Shreve \cite{Bertsekas_1978} for the discrete time case,
	where the main idea is to interpret control as the probability measures on the state space which lie in a topological space.
	The measurable selection theory justifies first the measurability of the value function,
	and further permits to ``paste'' a class of local controls into a global control by composition of probabiliy measures.

	For continuous time stochastic control problem, the DPP, in essential, is similarly based on the stability of the control under ``conditioning'' and ``pasting''.
	However, it becomes much more technical for a proof.
	To avoid technical questions, such as the measurability of the value function, one usually imposes sufficient conditions on the control problem which guarantees the continuity or semi-continuity of the value function,
	and then uses a separability argument to paste the controls, see e.g. Fleming and Soner \cite{Fleming_Soner} among many related works.
	For the same purpose, Bouchard and Touzi \cite{Bouchard_2011} proposed a weak dynamic programming by considering the semi-continuity envelop of the value function.
	In the same continuous time control context, El Karoui, Huu Nguyen and Jeanblanc \cite{ElKaroui_1987} interpreted the controls as probability measures on the canonical space of continuous paths and then proposed a general framework to derive the DPP using the measurable selection theorem.
	In the same spirit, Nutz and van Handel \cite{NutzHandel} proposed a general construction of time-consistence sublinear expectations on the canonical space, where their time-consistence property is a reformulation of the DPP.
	Using the same techniques, the DPP has been proved for some specific control problems, e.g. in Tan and Touzi \cite{TanTouzi}, Neufeld and Nutz \cite{NeufeldNutz} etc.
	
	\vspace{2mm}

	The first objective of this paper is to derive a general measurable selection theorem, following the approach developed in Dellacherie \cite{Dellacherie_1972} and in Dellacherie and Meyer \cite{Dellacherie_Meyer_1978}.
	The main idea is first to give an explicit construction of the selection in the ``compact'' case;
	then to extend it to the ``Borel'' case by approximation using capacity theory.
	Finally, it follows by an ``isomorphic'' argument that one obtains a general selection theorem.

	For a second objective, we propose a general framework for the DPP following \cite{ElKaroui_1987} and \cite{NutzHandel}.
	We consider the canonical space of c\`adl\`ag trajectories as well as its extension spaces,
	for a family of nonlinear operators, indexed by stopping times, on the functional space of the canonical space,
	We derive a time consistence property of the operators under appropriate conditions.
	When the operators are interpreted as control/stopping problems, the time consistence property turns to be the dynamic programming principle.
	In particular, this framework can be considered as a continuous time extension of the discrete time dynamic programming of Bertsekas and Shreve \cite{Bertsekas_1978}, or the gambling house model studied in Dellacherie \cite{Dellacherie_MaisonJeux}, Maitra and Sudderth \cite{MaitraSudderth}.
	In our accompanying paper \cite{ElKaroui_Tan2}, we shall show that this framework is convient to study general stochastic control problems.

	\vspace{2mm}

	The rest of the paper is organized as follows.

	In Section \ref{sec:capacity}, we follow Dellacherie \cite{Dellacherie_1972} to derive a measurable selection theorem.
	We first give in Section \ref{subsec:capacity} a brief introduction to the capacities theory, including Choquet's capacitability theorem and in particular a projection capacity which induces a measurable projection result.
	Next in Section \ref{subsec:MeasurSelec}, we construct a measurable selection using the ``debut'' of the set in a product space.
	Then by an ``isomorphism'' argument, we obtain immediately a general measurable selection theorem.
	We also cite Jankov-von Neunman's analytic selection theorem in the end of the section.

	In Sectione \ref{sec:DPP_discrete}, we first recall some facts on probability kernel as well as the composition and disintegration of probability measures,
	then we present the classical approach to derive the dynamic programming principle for discrete time control problems using measurable selection.

	In Section \ref{sec:DPP}, we introduce an abstract framework for the dynamic programming principle (DPP).
	A detailed discussion on the canonical space of c\`adl\`ag trajectories is provided in Section \ref{subsec:canonical_space}.
	Then in Section \ref{subsec:TC}, we introduce a family of nonlinear operators defined on the functional space of the canonical space, and derive a time consistency property using measurable selection theorem.
	We also introduce an enlarged canonical space motivated by optimal stopping problems.
	In particular, the time consistency property turns to be the DPP when the nonlinear operators are induced by control rules as shown in Section \ref{subsec:abstr_dpp}.

	\vspace{2mm}

\noindent {\bf Notations}.
  We provide here some frequently used notations.\\
  \rmi The notation $E$ refers in general to a metric space, but usually $E$ is assumed to be a Polish space, i.e. a complete (every Cauchy sequence has a limit in itself) and separable ($E$ has a countable dense subset) metric space. Sometimes, $E$ is also assumed to be a locally compact Polish space.
 \begin{itemize}
        \item[$-$] $\Kc(E)$ (or simply by $\Kc$) denotes the class of all compact subsets of $E$, i.e.
	    	\be \label{eq:Kc}
            \Kc ~=~ \Kc(E) ~:=~ \{~ \mbox{All compact subsets of}~ E ~\}.
        \ee
 \item[$-$]  $\Bc(E)$ denotes the Borel $\sigma$-field generated by all the open sets (or closed sets) in $E$.
        \item[$-$] $\Mc(E)$ (resp. $\Pc(E)$) denotes the collection of finite positive measures (resp. probability measures) on $(E, \Bc(E))$.
        When $E$ is a Polish space, $\Mc(E)$ and $\Pc(E)$ are both Polish spaces equipped with the topology of the weak convergence, i.e. the coarsest topology making
			\b*
				\mu \in \Pc(E) &\mapsto& \mu(\varphi) ~:=~ \int_E \varphi(x) \mu(dx) \in \R
			\e*
			continuous for every bounded continuous function $\varphi \in C_b(E)$. Moreover, the Borel $\sigma$-field of the Polish space $\Pc(E)$ is generated by the maps  $\mu \mapsto \mu(\varphi)$ with all $\varphi \in C_b(E)$.
  \end{itemize}
\rmii Abstract spaces are referred by different notations:
 \begin{itemize}
   \item[$-$] for instance we  use $X$ to refer to an abstract set, where $2^X$ denotes the collection of all subsets of $X$. Let $\Xc$ be a $\sigma$-field on $X$, then $(X,\Xc)$ is called a measurable space, and $\Lc(\Xc)$ denotes the collection of all $\Xc$-measurable functions.
    \item[$-$] When the context is more probabilistic as in  Section \ref{sec:capacity}, we use the classical notation $(\Om, \Fc)$ in place of  $(X,\Xc)$.
   \item[$-$] Nevertheless, in Section \ref{sec:DPP} without ambiguity, $\Om := D(\R^+, E)$ denotes the canonical space of all $E$-valued c\`adl\`ag paths on $\R^+$ with a Polish space $E$.
     \end{itemize}
\rmiii We use the usual convention that the supremum (resp. infimum) over an empty subset is $-\infty$ (resp. $\infty$),  i.e. $\sup \emptyset = -\infty$ (reps. $\inf \emptyset = \infty$).







\section{Capacity theory and measurable selection}
\label{sec:capacity}

	The goal of this section is to give a simple and straightforward proof of the projection and the selection theorems using the capacity theory. The most remarkable phenomena is the role of the negligible or measure zero sets, which allow to a very general measurability result up to a negligible set.
	
	To illustrate the idea, let us consider the product space $[0,1] \x [0,1]$ and a Borel subset $A \subseteq [0,1] \x [0,1]$.
	The question is whether the projection set
	$$\pi_{[0,1]}(A) ~:=~ \big\{x \in [0,1] ~: \exists y \in [0,1], ~\mbox{s.t.}~ (x,y) \in A \big\}$$
	is measurable in [0,1].
	When $A$ is compact, $\pi_{[0,1]}(A)$ is still compact and hence a Borel set.
	When $A$ is only Borel, we can approximate it from interior by a sequence of compact subsets, using Choquet's capacitability theorem, and then deduce a measurability result up to a negligible set.

	In the following of this section, we start with a brief introduction to the capacity theory up to  Choquet's capacitability theorem.
	The capacity is formulated first in a topological approach, and then in an abstract way.
	By considering a projection capacity, we deduce naturally a selection measurable theorem in Section \ref{subsec:MeasurSelec}.

\subsection{Paving and abstract Choquet's capacity}
\label{subsec:capacity}
\subsubsection{Capacity on topological space}
Let us start by a simple example of capacity on a metric space $E$ to illustrate the differences between measures and capacities, and to motivate the future abstract development.\\[-8mm]
\paragraph{\small Outer measure and outer/inner regularity}
\label{oiregularity}

	The simplest example of capacity is the outer measure on a metric space $E$ equipped with its Borel  $\sigma$-field $\Bc(E)$.
  Let $\mu$ be a finite positive measure on $(E, \Bc(E))$, then its outer measure, acting on subsets of $E$, is defined by
	\b*
		\forall A \subseteq E, ~~~ \mu^*(A) &:=& \inf ~\big \{ \mu(C) ~:  C \text{ open set containing } A  \big\}.
	\e*
	Since the  outer  measure $\mu^*$ coincides with the  measure $\mu $ on the Borel $\sigma$-field,
	it follows the {\em outer regularity by open sets} of the measure $\mu(A)$ for all $A \in \Bc(E)$.
	On the other hand, thanks to the topological structure of $E$,
	it is well-known that an approximation by below by {\em closed sets} also holds true, i.e.
    \b*
        \forall B \in \Bc(E), && \mu(B) ~=~ \sup ~\big\{ \mu(C) ~: C ~\mbox{is closed, contained in}~ B \big\}.
    \e*
	When in addition $E$ is a Polish space, the inner approximation by closed sets can be replaced by the compact sets,
	known as {\em $\Kc$-inner regularity}, i.e.
	\be \label{eq:inner_regular}
	   \forall B\in \Bc(E),&& \mu(B) ~=~ \sup ~\big \{ \mu(K) ~: K \in \Kc(E) \text{ contained in } B \big \}.
    \ee
  In 1955, Choquet \cite{Choquet_1955} observed that the inner regularity property does not depend on the additivity of $\mu$ on $E$, but only on its monotonicity and sequential continuity.
  These facts yield to the notion of Choquet's capacity:
  \begin{Definition}
    A capacity $C$ on Polish space $E$ is a mapping $C: 2^E \to \R^+$ which is monotone (i.e. $A \subseteq B \Rightarrow C(A) \le C(B)$),
    continuously increasing (i.e. $A_n\uparrow A \Rightarrow C(A_n) \uparrow  C(A)$),
    and continuously decreasing on the compact sets, that is $K_n\in \Kc,\, K_n\downarrow K \Rightarrow C(K_n)\downarrow  C(K)$.
  \end{Definition}

    Besides the outer measures, the supremum of the out measures over a compact set of measures is a typical example of capacity.
	\begin{Proposition}[Supremum of outer-measures] \label{supmeasure}
    Let $E$ be a Polish space with compact class $\Kc$.
    Denote by $\Pc(E)$ the collection of all probability measures on $(E,\Bc(E))$ and by $\Pc_{K}(E)$ a compact subset of $\Pc(E)$ for the weak topology.
    Then, the  set function $\I$ defined as the supremum of the family of outer measures is a Choquet's capacity defined below \eqref{eq:inner_regular},
		\b*
			\forall A \subseteq E,~~~ \I(A) :=\,\sup_{\mu \in \Pc_{K}(E)} \mu^*(A).
		\e*
	\end{Proposition}
	\proof First, it is obvious that $\I$ defined above takes values in $\R^+$ and is monotone.
	Moreover, as the supremum of a family of outer measures, $\I$ has clearly the increasing continuity.
	To deduce the decreasing continuity of $\I$ on $\Kc$, we use the compactness property of $\Pc_K(E)$.
        Suppose that $(K_n)_{n \ge 1}$ is a decreasing sequence of nonempty compact sets of $E$, then $K := \cap_n K_n$ is still nonempty and compact.
        By extracting subsequence if necessary, there exists a sequence $(\mu_k, K_{n_k})_{k \ge 1}$ in $\Pc_{K}(E) \x \Kc$ such that $\mu_k \to \mu$,
        and $\mu_k(K_{n_k}) \searrow \inf_n \I(K_n)$ as $k \to \infty$. Since $K_{n_k}$ are all closed sets,
        it follows that $\mu(K_{n_{k_0}}) \ge \limsup_{ k \to \infty} \mu_k(K_{n_{k_0}}) \ge \lim_{ k \to \infty} \mu_k(K_{n_k}) = \inf_n \I(K_n)$ for every $k_0 \ge 1$.
        Hence $\I(K) \ge \mu(K) \ge \inf_n \I(K_n)$, and with the other obvious inequality, we conclude the proof for compact $\Pc_{K}(E)$. \qed

    \begin{Remark}
        {\rm Choquet's theorem proves the $\Kc$-inner regularity of any topological capacity on the Borel sets, that is,
        $\forall B \in \Bc(E), ~C(B ) = \sup \big \{ C(K) ~: K \in\Kc,  \,K \subseteq B \big \}$.
        Hence, $\I$ defined above has the $\Kc$-inner regularity property.
        In fact, as invoked in Denis, Hu and Peng \cite{DenisHuPeng} (see e.g. their Theorem 2), the supremum of a family of capacities has always $\Kc$-inner regularity.}
    \end{Remark}

	\begin{Remark} \label{rem:capacity_topo}
		{\rm
		\rmi These kind of sublinear capacities has always been a major tool not only in potential theory with Choquet's work, but also in many areas as stochastic processes, optimization and games theories.
        A recent renewed interest comes from the theory of non-linear expectation on canonical path space (see Sections 3 and 4 and the reference therein). \\
		\rmii A functional version of capacity $C$ is defined in Dellacherie \cite{Dellacherie_1972b} as a mapping from the space of all functions on a Polish space $E$ to $\R$,
        which is monotone (i.e. $f \le g \Rightarrow C(f) \le C(g)$), continuously increasing on general functions (i.e. $f_n \uparrow f \Rightarrow C(f_n) \uparrow C(f)$)
        and continuously decreasing on upper semicontinuous functions,
        that is for decreasing sequence of upper semicontinuous functions $g_n \downarrow g$, one have $ C(g_n) \downarrow C(g)$.
		}
	\end{Remark}


\subsubsection{Paving and Abstract capacity}
\label{subsubsec:paving}
	In 1959, Choquet \cite{Choquet_1959} generalized the above inner regularity \eqref{eq:inner_regular} into an abstract form, named ``capacitability'',
    with the capacity theory, where in general the role of compact sets is played by a family of subsets with the same stability properties as $\Kc.$
    Here, we follow Dellacherie \cite{Dellacherie_1972} to define the paving (``pavage'' in French), which differs slightly from the definition in Dellacherie and Meyer \cite{Dellacherie_Meyer_1978}.

	\begin{Definition}[\bf Paving as abstract version of $\Kc$] \label{def: paving}Let $X$ be a set.\\
    \rmi A {\bf paving} $\Jc$ on $X$ is a collection of subsets of $X$ which contains the empty set $\emptyset$ and which is stable under finite intersections and finite unions. The couple $(X,\Jc)$ is called a {\bf paved space}.\\
		\rmii A  paving stable under countable intersections is called a {\bf $\delta$-paving}; and a $\delta$-paving stable under countable unions is called a {\bf mosaic}.\\
		\rmiii As usual, given a class $\Ic$ of subsets of $X$, $\Ic_p$ (resp. $\Ic_{\delta}$, $\widehat{\Ic}$, $\sigma(\Ic)$) denote the smallest paving (resp. $\delta$-paving, mosaic, $\sigma$-field) containing  $\Ic$. \\[-6mm]
	\end{Definition}

	Two typical examples of pavings are the  class of compact sets $\Kc$ in a metric space $E$, and the $\sigma$-field $\Fc$ in an abstract space $\Omega$.
	In particular, they are both $\delta$-pavings, i.e. $\Kc_\delta=\Kc, \>\Fc_\delta=\widehat\Fc=\Fc$.
	Now, with a paved space, we can extend the definition of capacity to an abstract framework. 	
	\begin{Definition}\label{def:capacity}
		Let $(X, \Jc)$ be a paved space,
		a $\Jc$-capacity $\I$ is a mapping from $2^X$ to $\Rb := \R \cup \{ \infty \}$ satisfying:\\
		\rmi  $\I$ is increasing, i.e. if $A \subseteq B$, then $ {\I}(A) \le {\I}(B)$;\\
        \rmii $\I$ is sequentially increasing on $2^X$, i.e. for every increasing sequence $(A_n)_{n \ge 1}$ in $2^X$, $\>\I( \cup_n A_n ) = \sup_n \I(A_n)$;\\
		\rmiii $\I$ is sequentially decreasing on $\Jc$, that is for every decreasing sequence $(A_n)_{n \ge 1}$ in $ \Jc$, $\>\I ( \cap_n A_n) = \inf_n \I(A_n)$.
	\end{Definition}
	We provide directly the main theorem on the capacitability, due to Choquet \cite{Choquet_1959}; and we refer to Theorem I.31 of Dellacherie \cite{Dellacherie_1972} for a detailed proof.

	\begin{Theorem}\label{theo:Choquet}(Choquet)
		Let $(X, \Jc)$ be a paved space and $\I$ be a $\Jc$-capacity. Then any element $C$ of  the mosaic $\widehat{\Jc}$ is  $\Jc$-capacitable, i.e. its capacity can be approximated from below by elements of $\Jc_{\delta}$,
		\be \label{eq:capacitability}
			\forall C \in \widehat{\Jc}, ~~~~ \I(C) ~=~ \sup\{ ~ \I(K)~ : K \in \Jc_{\delta}, ~K \subseteq C ~\}.
		\ee
	\end{Theorem}

	\begin{Remark}
		{\rm
			\rmi Let $E$ be a metric space with its compact paving $\Kc$,
			then every topological capacity defined above Proposition \ref{supmeasure} is clearly a $\Kc$-capacity.\\
			\rmii Let $(X, \Xc)$ be an abstract measurable space, $\mu$ be a finite positive measure, the outer measure $\mu^*$ defined by\\
			\centerline{$\forall A \in 2^X,~~~ \mu^*(A) ~:=~ \inf ~\{ \mu(B) ~: A \subseteq B ~~\mbox{and}~~ B \in \Xc \}$~~~~~~~~~~}\\
			is then a $\Xc$-capacity, since $\mu^*$ coincides with $\mu$ on the $\sigma$-field $\Xc$.\\
			\rmiii For the functional capacity given in Remark \ref{rem:capacity_topo},
            the capacitability property \eqref{eq:capacitability} turns to an approximation of a Borel function by an increasing sequence of upper semicontinuous functions as shown in Dellacherie \cite{Dellacherie_1972b}.
		}
	\end{Remark}

\subsubsection{  Product paving and projection capacity}

	Let us turn back to the example given in the beginning of the section, the question is whether the projection set $\pi_{[0,1]}(A) \subseteq [0,1]$ is measurable given a Borel set $A \subset [0,1] \x [0,1]$.
	In fact, Lebesgue has argued that $\pi_{[0,1]}(A)$ is obvious a Borel set; and it is Souslin who pointed out that $\pi_{[0,1]}(A)$ may not be Borel,
    who was then motivated to develop the theory of analytic set (see also Section \ref{subsubsec:MeasSelec_analytic}).

	Using the capacity theory, one can partially reply to the question, that is, $\pi_{[0,1]}(A)$ is measurable w.r.t. the completed $\sigma$-field.
	To this end, we shall consider a product paving as well as a projection capacity.\\[-8mm]

\paragraph{\small Product paving and projection operator}\label{productpaving}

	Let $(\Om, \Fc)$ be a measurable space, $E$ be a {\bf locally compact} Polish space equipped with its compact paving $\Kc$.
    The product paving $\Ic_p=\Fc \otimes_p \Kc$ is generated by the class $\Ic := \Fc \oplus \Kc$  of all rectangles $B\times K$ with  $B\in \Fc\> \mbox{and} \>K\in \Kc$.
    Since $\Ic $ is stable by finite intersections,
    any element in the paving $\Ic_p$ is a finite union of elements in $\Ic$, and the $\delta$-paving $\Ic_\delta =\Fc \otimes_\delta \Kc$
    is composed by all elements which are the decreasing countable intersection of decreasing elements of $\Ic_p$.
    Moreover, thanks to the separability and local compactness of $E$, for every $A = B \x K \in \Ic$ such that $B \in \Fc$ and $K \in \Kc$,
    the complement $A^c = (A^c \x E ) \cup (\O \x K^c)$ lies in the mosaic $\widehat{\Ic} = \Fc \widehat{\otimes}_p \Kc$ generated by $\Ic= \Fc \oplus \Kc$.
    It follows that $\widehat{\Ic} = \Fc \widehat{\otimes}_p \Kc$ coincides with the product $\sigma$-field $\sigma(\Ic) = \Fc \otimes \sigma(\Kc)= \Fc \otimes \Bc(E)$.
%
   \begin{Proposition}\label{prop:projection} Let us consider the product space $\Omega \times E$ equipped with the paving $\Ic_p$  generated by the rectangles
   $ \Fc\oplus\Kc$.
	The  projection operator $\pi_{\O}$ maps subsets of $\Omega \times E$ to subsets of $\Omega$,
    \be \label{proj}
        \pi_{\O}(A) ~:=~ \big\{~ \om ~: \exists ~x\in E ~\text{such that }\>(\om,x) \in A ~\big \}.
    \ee
        \rma $\pi_{\O}$ is continuously increasing, i.e. \\
        \centerline{$A \subset B ~\Rightarrow~ \pi_{\O}(A)  \subset \pi_{\O}(B) $ ~~and~~ $A_n\uparrow A ~\Rightarrow~ \bigcup_n \pi_{\O}(A_n)= \pi_{\O}(A)$.~~~}\\
        \rmb  $\pi_{\O}$  is continuously decreasing on $\Ic_\delta$, i.e. $L_n \in \Ic_\delta, \, L_n\downarrow L \Rightarrow \bigcap_n\pi_{\O}(L_n)=  \pi_{\O}(L)$.\\
        \rmc The class $ \pi_{\O}(\Ic_\delta)$ is a $\delta$-paving on $\Omega$.
   \end{Proposition}
   \proof  It is enough to prove assertion \rmb since the other two items are obvious.
        Let $L_n \in \Ic_\delta$ be a decreasing sequence, such that $(L_n)_{n \ge 1}$ are all nonempty.
        First, suppose that $\om \in \pi_{\Om}(L)$, then there is $x \in E$ such that $(\om, x) \in L \subset L_n$ for all $n \ge 1$.
        Therefore, $\om \in \pi_{\Om}(L_n)$ for all $n \ge 1$, and hence $\om \in \bigcap_n \pi_{\Om}(L_n)$.
        On the other hand, suppose that $\om \in \bigcap_n \pi_{\Om}(L_n)$.
        Since all elements in $\Ic_{\delta}$ can be represented as the decreasing countable intersections of decreasing elements of $\Ic_p$,
        it follows that the section set $L_n^{\om} := \{ x ~:(\om, x) \in L_n \}$ is nonempty and compact.
        Therefore, $\bigcap_n L_n ^{\om}$ is nonempty, which implies that $\om \in \pi_{\O}( \bigcap_n L_n) = \pi_{\Om}(L)$.
	\qed\\[-8mm]

\paragraph{Projection capacity} \label{projectioncapacity}

    We now introduce the projection capacity on the product paving defined above, which plays an essential role to derive a measurable selection theorem.

	Let $(\O, \Fc, \P)$ be a probability space, $\P^*$ be the outer measure, $E$ be a locally compact Polish space with its compact paving $\Kc$.
	The product space $\Omega \times E$ is equipped with the paving $\Ic_p=\Fc \otimes_p \Kc$  generated by all the rectangles in $\Fc\times\Kc$.
    Using the projection operator $\pi_{\Om}$ introduced in \eqref{proj}, we define a set function on the subsets of $\O \x E$ by
	\be \label{eq:capacity_I}
		\forall A \subseteq \O \x E, ~~~ \I_\Omega(A) &:=& \P^*( \pi_{\O}(A) ).
	\ee

	\begin{Theorem}\label{theo:measurable_projection}
    \rmi The set function $\I_\Omega$ is a $\Fc \otimes_p \Kc$-capacity, we call it the projection capacity.\\
    \rmii For every $A \in   \Fc \otimes \Bc(E)$, there exists an increasing sequence $(L_n)_{n \ge 1}$ in $(\Fc \otimes_p \Kc)_{\delta}$
so that $\pi_{\O}(L_n) \in \Fc$ and
		\be \label{eq: approx cap}
			L_n \subseteq A,\quad \P(\pi_{\O}(L_n)) ~~\le~~ \P^*( \pi_{\O}(A) ) ~~\le~~  \P(\pi_{\O}(L_n)) ~+~ \frac{1}{n}.
		\ee
    \rmiii It follows that $\P^*( \pi_{\O}(A) )=\P(\cup_n \pi_{\O}(L_n))$, where $\pi_{\Om}(A) \supseteq \cup_n \pi_{\O}(L_n) \in \Fc$.
    In other words, the projection $\pi_{\O}(A)$ differs from the measurable set $\cup_n \pi_{\O}(L_n) $ by a $\P$-negligible set, and then is $\Fc$-measurable if the $\sigma$-field $\Fc$ is $\P$-complete.
\end{Theorem}
\proof  \rmi There is no difficulty in proving that $\I_\Omega$ is a capacity, given the properties of the projection operator $\pi_\Omega$ stated in Proposition \ref{prop:projection} and the fact that $\P^*$ is a capacity.
        \rmii Since the mosaic generated by the paving $\Fc \otimes_p \Kc$ is the product $\sigma$-field $\Fc \otimes \Bc(E)$, by Choquet's theorem (Theorem \ref{theo:Choquet}),
        we can construct a sequence $(L_n)_{n \ge 1}$ in $(\Fc \otimes_p \Kc)_{\delta}$ such that \eqref{eq: approx cap} holds true.
        Further, by simple manipulation, the inner approximating sequence $(L_n)_{n \ge 1}$  may be chosen non decreasing.
        \rmiii The last item is an immediate consequence of (ii).
        \qed

    \vspace{2mm}

	The functional version of the above projection result is related to the measurability of the supremum of a family of random variables, where the index parameter lies in a topological space.

	\begin{Corollary} \label{Coro:Measurability_Sup_K}
	Let $(\O, \Fc,\P)$ be a probability space and $E$ be a  locally compact Polish space.  We consider  a measurable function $ f(\om , x)$  defined on $ ( \Omega \times E, \Fc\otimes \Bc(E) )$, and $A \in \Fc\otimes \Bc(E)$.
    Then the function $g$ defined by
		\be \label{eq:supdef}
			g(\om) &:=& \sup \{~ f(\om, x) ~: (\om, x) \in A ~\}, ~\mbox{with the convention}~ \sup \emptyset = -\infty,
		\ee
		is $\Fc$-measurable whenever $\Fc$ is $\P$-complete.
	\end{Corollary}
	
	\proof  Given a constant $c \in \R$, let us define $B_c := \big \{(\om, x) \in A ~: f(\om, x) > c \big \}$, and  $C_c := \{ \om ~: g(\om) > c \}.$
	Clearly, $B_c \in \Bc(E) \otimes \Fc$ and $C_c = \pi_{\O}(B_c)$. When $\Fc$ is $\P$-complete, it follows from Theorem \ref{theo:measurable_projection} that $C_c \in \Fc$, and hence $g$ is $\Fc$-measurable.
	\qed

	\begin{Remark}{\rm
        \rmi To see that Corollary \ref{Coro:Measurability_Sup_K} is a functional version of the measurable projection theorem (Theorem \ref{theo:measurable_projection}),
        we can consider the function $f(\om,x) := \1_{A}(\om,x)$ for $A \in \Fc \otimes \Bc(E)$,
        then Corollary \ref{Coro:Measurability_Sup_K} implies that $g(\om) = \1_{\pi_{\O}(A)}(\om) - \infty \1_{\pi_{\O}(A)^c}(\om)$ is measurable, which induces the measurability of $\pi_{\O}(A)$.\\
		\rmii We notice that in Theorem \ref{theo:measurable_projection} and Corollary \ref{Coro:Measurability_Sup_K},
        although the approximation sets $(L_n)_{n \ge 1}$ depend on the probability $\P$ on $(\O,\Fc)$, the projection $\pi_{\O}(A)$ and the supremum function $g$ are independent of $\P$.
        Therefore, the completeness condition of $(\O, \Fc, \P)$ can be relaxed.
        We shall address this issue later in Section \ref{subsubsec:Meas_Sele_univ}.\\
        \rmiii When $(\Om, \Fc, \P)$ is a complete probability space, the supremum function $g$ is $\Fc$-measurable.
        However, $g$ may not be the essential supremum of the family $(f(\om,x))_{x \in E}$ in general.
        For example, let $\Om = E = [0,1]$, $A = [0,1] \x [0,1]$ and $f(\om,x) = 1_{\om = x}$, it follows that $g(\om) \equiv 1$ and the essential supremum of the family $(f(\om,x))_{x \in E}$ under the Lebesgue measure is $0$.
        }
	\end{Remark}

\subsection{Measurable  selection theorem}
\label{subsec:MeasurSelec}
    Now, we adopt the presentation of Dellacherie \cite{Dellacherie_1972} to show how to deduce the measurable selection theorem using the projection capacity,
    where the case of $\R^+$ is considered in a first step and an extension is obtained by the isomorphism argument.
    In the end, we also cite another measurable selection theorem presented in Berstekas and Shreve \cite{Bertsekas_1978}, or in Bogachev  \cite{Bogachev}.

\subsubsection{Measurable selection theorem: the $\R^+$ case}
\label{subsubsec:MeasurSelec_R}

	The order structure of $\R^+$ allows us to deduce very easily a measurable selection result from the measurable projection theorem (Theorem \ref{theo:measurable_projection}).
    It is a key result  in the study of stochastic processes indexed by continuous time that has motivated the  following presentation by Dellacherie \cite{Dellacherie_1972},
    with the familiar notion of the debut of the sets in the product space $\Omega\times  \R^+$.
    We shall start with a {\bf complete} probability space.

	\begin{Proposition} \label{prop:measurability_debut}
		Let $(\O, \Fc, \P)$ be a complete probability space and $A$ a subset of the product space $\Omega \x \R^+$, measurable with respect to  the product $\sigma$-field $ \Fc \otimes \Bc(\R^+) $. \\
		\rmi The debut  $D_A$ of $A$ is  defined by
		\be \label{eq:debut_A}
			D_A(\om) &:=& \inf \big \{ t \geq 0 ~: (\om,t) \in A \big\},
		\ee
		which yields to the identity $\{\om ~: D_A(\om)< \infty \}=\pi_{\O}(A)$.
		Then the debut $D_A$ is a $\Fc$-measurable random variable.\\
		\rmii In addition, let $\F = (\Fc_t)_{t \ge 0}$ be a filtration on $\O$ satisfying the usual conditions (i.e. $\Fc_0$ contains all $\P$-null sets and $\F$ is right continuous),
        and $A$ a $\F$-progressive set in $ \O \x \R^+$ (i.e. $\forall t \ge 0,\> A \cap (\O \x [0,t] ) \in  \Fc_t \otimes \Bc([0,t])$), then $D_A$ is a $\F$-stopping time.
	\end{Proposition}
	\proof \rmi The measurability of $D_A$ follows immediately by Theorem \ref{theo:measurable_projection} as well as the fact that for every $ t \in \R^+$, $\big\{ \om ~: D_A(\om) < t \big \} ~=~ \pi_{\O} \big( A \cap (\O \x [0,t)) \big)$.\\
		\rmii When $A$ is progressively measurable, it follows by similar arguments that for every $t \in \R^+$, $ \{ \om ~: D_A(\om) < t \} \in \Fc_t$.
        Since the filtration $\F$ satisfies the usual conditions, it follows that $D_A$ is a $\F$-stopping time. \qed

    \vspace{2mm}

	Let us stay in the context of Proposition \ref{prop:measurability_debut} and define the graph $[[D_A]]$ of $D_A$ as a subset of $\O \times \R^+ $ by
	  \be \label{eq: graph}
			[[D_A]] &:=& \big \{( \om,D_A(\om)) \in \O \times \R^+  \ \big \} ~=~ \big \{( \om,D_A(\om)) ~: D_A(\om) <+ \infty \big \}.
	  \ee
	If $[[D_A]] \subseteq A$, which is the case when $A^{\om}$ is compact for every $\om$, then it is clear that the graph of debut $D_A$ gives a measurable section of set $A$.
    Otherwise, we can easily overcome this difficulty using the approximation techniques given by Choquet's theorem (Theorem \ref{theo:Choquet}), and hence establish the following measurable selection theorem.

	\begin{Theorem}\label{theo:measurable_section_t}
		Let $(\O, \Fc, \P)$ be a complete probability space, $A \in \Fc \otimes \Bc(\R^+)$ be a measurable subset in $\O \x \R^+$.
        Then there exists a $\Fc$-measurable random variable  $T$ taking values in $[0,\infty]$ such that
		\be \label{eq:measruable_section_t}
			[[T]] = \big\{ \big( \om, T(\om) \big) \in \O \x \R^+ \big\} \subseteq A ~~\text{and}~~\big\{ \om ~: T(\om) < \infty \big \} ~=~ \pi_{\O}(A).
		\ee
	\end{Theorem}
	\proof Let $\Fc \otimes_p \Kc$ be the product paving defined in the end of Section \ref{subsubsec:paving}, then by the measurable projection theorem (Theorem \ref{theo:measurable_projection}),
        there exists an increasing sequence $(L_n)_{n \ge 1}$ in $(\Fc \otimes_p \Kc)_{\delta}$ such that $ L_n \subseteq A ~~\text{and}~ \I_{\Om} (A) ~\le~ \I_{\Om} (L_n) + \frac{1}{n}$,
        where $\I_{\Om}$ is the projection capacity defined in \eqref{eq:capacity_I}.
        The debut $D_{L_n}$ of this set is $\Fc$-measurable by Proposition \ref{prop:measurability_debut}.
        Moreover, $[[D_{L_n})]] \subseteq L_n \subseteq A$ since the section set $L_n^{\om}$ is compact for every $\om \in \pi_{\O}(L_n)$.
        Let $T_1 := D_{L_1}$, $T_{n+1} := T_n \1_{T_n < \infty} + D_{L_{n+1}} \1_{T_n = \infty} $ and $T_{\infty} := \lim_{n \to \infty} T_n$.
        Then clearly, $[[T_{\infty}]] \subseteq A$ and $T_{\infty}(\om) < \infty$ for every $\om \in \cup_{n=1}^{\infty} \pi_{\Om}(L_n) \subseteq \pi_{\Om} (A)$.
        Using the fact that $(\O, \Fc, \P)$ is complete and $\pi_{\Om}(A) \setminus (\cup_{n=1}^{\infty} \pi_{\Om}(L_n))$ is $\P$-negligible,
        it is easy to construct the required $T$ from $T_{\infty}$. \qed

	\vspace{2mm}

	In the general theory of processes, with a complete probability space $(\O, \Fc, \P)$  equipped with a filtration  $\F =(\Fc_t)_{t \ge 0},$ satisfying the usual conditions,
    similar results are shown for the paving and the optional (resp. predictable) $\sigma$-field on the space $\Omega \x \R^+$ generated by the subsets $[[S,\infty)) := \{(\om, t) ~: t \ge S(\om) \}$,
    where $S$ are stopping times (resp. predictable stopping times).
    In this case, for an optional (resp. predictable) set $A$, one has an approximation sequence $(L_n)_{n \ge 0}$ whose debuts are $\P$-a.s. stopping times (resp. predictable stopping times),
    and the section theorem becomes: for any optional (resp. predictable) set $A$, and $\eps >0$,
    there exists a stopping time (resp. predictable stopping time) $T$ such that $[[T]] \subseteq A$ and $\P(T<\infty)\geq \P(\pi_{\O}(A))-\eps$.
    An immediate consequence is the following result:

	\begin{Proposition}
		Let $X$ and $Y$ be two optional (resp. predictable) nonnegative processes such that for any stopping times $T$ (resp. predictable stopping times), $\E[X_T \1_{T<\infty}]=\E[Y_T \1_{T<\infty}]$.
        Then X and Y are indistinguishable, i.e. $\P(X_t=Y_t, \forall t\in \R^+)=1$;
        in other words the projection set $\pi_{\Om} (A)$ of $A := \{(\o,t) ~: X_t(\o)\neq Y_t(\o)\}$ is a  $\P$-negligible set.
	\end{Proposition}
	\proof It is enough to show that for every $\eps$, the optional set $A_{\eps} = \{(\omega,t) ~: X_t(\om) >  Y_t(\om) + \eps \}$ satisfies $\P(\pi_{\Om}(A_{\eps})) = 0$.
    It is true since otherwise, there is a measurable section $T_{\eps}$ of $A_{\eps}$ such that $\E [ X_{T_{\eps}} 1_{T_{\eps} < \infty} ] \ge \E [ Y_{_{\eps}} 1_{T_{\eps} < \infty} ] + \eps \P(\pi_{\Om}(A_{\eps}))$,
    which contradicts the assumption of the proposition. \qed

\subsubsection{Universally measurable selection}
\label{subsubsec:Meas_Sele_univ}

	It is clear that when $(\O, \Fc, \P)$ is not complete, it can be completed by its outer measure $\P^*$.
    Let us denote by $\Nc^{\P}$ the collection of all subsets $N \subseteq \O$ such that $\P^*(N) = 0$, the elements in $\Nc^{\P}$ is called $\P$-negligible set.
    Then $\Fc^{\P} = \Fc \bigvee \Nc^{\P} := \sigma( \Fc \cup \Nc^{\P} )$ is called the completed $\sigma$-field with respect to probability $\P$.
    The probability measure $\P$ is extended uniquely on $\Fc^{\P}$ by its outer measure $\P^*$, which is also a probability measure on $(\O, \Fc^{\P})$.
    Let us denote the new probability space by $(\O, \Fc^{\P}, \P)$.

	Any $\P$-negligible set in $\Fc^{\P}$ is contained in a $\P$-negligible set in $\Fc$.
    Moreover, let $A \in \Fc^{\P}$, there are $B, N \in \Fc$ such that $\P(N) = 0$ and $A = B \Delta N := (B\setminus N) \cup (N \setminus B)$.
    More generally, given a $\Fc^{\P}$-measurable random variable $S : \O \to \R$,
    we can construct a $\Fc$-measurable random variable $\tilde{S}$ by the approximations with step functions such that $S = \tilde{S}$, $\P$-almost surely.
    Then the measurable selection $T$ in Theorem \ref{theo:measurable_section_t} can be chosen in a ``almost surely'' sense in order to make $T$ measurable w.r.t. $\Fc$, even if the probability space $(\O, \Fc, \P)$ is not complete.
	
	Finally, given a measurable space $(\O, \Fc)$ and $A \in \Fc \otimes \Bc(\R^+)$, we can always complete $\Fc$ by $\Fc^{\P}$ with a probability measure $\P$ and then get $\pi_{\O}(A) \in \Fc^{\P}$.
    Nevertheless, the definition of $\pi_{\O}(A)$ does not depend on any probability measure $\P$.
    It follows that it is $\Fc^{\P}$-measurable for every probability measure $\P$ on $(\O, \Fc)$, which implies that $\pi_{\O}(A) $ is universally measurable,
    i.e. measurable with respect to the universal completion of $\Fc$ defined as follows.
	\begin{Definition}
		Let $(\O, \Fc)$ be a measurable space, the universal completion of $\Fc$ is the $\sigma$-field defined as the intersection of $\Fc^{\P}$ for all probability measures $\P \in \Pc(\O)$ on $(\O, \Fc)$, i.e.
		\be
			\Fc^U &:=& \bigcap_{\P \in \Pc(\Om)} \Fc^{\P}.
		\ee
	\end{Definition}
	By the same arguments, it follows that the supremum function $g$ defined by \eqref{eq:supdef} in Corollary \ref{Coro:Measurability_Sup_K} is also $\Fc^U$-measurable.
    But the construction of the selection $T$ in Theorem \ref{theo:measurable_section_t} depends on the measure $\P$.
    Nevertheless, it is still possible to get a universally measurable selection by other methods.
    In fact, a more precise characterization is obtained as analytic measurable selection when $\Om$ is a topological space with its Borel $\sigma$-field (see also Section \ref{subsubsec:MeasSelec_analytic} below).
    Dellacherie and Meyer \cite{Dellacherie_Meyer_1978} gives further an extension in the abstract measurable context in their Theorem III-82,
    where the main idea is to reduce the general $\sigma$-field case to a separable $\sigma$-field case, which is further linked to the topological context of the analytic measurable selection theorem.
    We accept their result and give the following conclusion.
	
	\begin{Proposition} \label{prop:universal_section_t}
		Let $(\Om, \Fc)$ be a measurable space, $A \in \Fc \otimes \Bc(\R^+)$, then $\pi_{\Om}(A) \in \Fc^U$ and there is a $\Fc^U$-measurable selection $T$,
        i.e. $\{ \om ~:T(\om) < \infty \} = \pi_{\Om}(A)$ and $[[T]] := \{ (\om, T(\om)) \in \Om \x \R^+ \} \subseteq A$.
		Suppose, in addition, that $f : \Om \x \R^+ \to \R$ is a measurable and $g$ is the supremum function	
        $g(\om) := \sup \{ f(\om, x) ~: (\om, x) \in A \}$, then $g: \Om \to \R \cup \{ - \infty, + \infty \}$ is $\Fc^U$-measurable.
	\end{Proposition}
	

\subsubsection{Measurable selection theorem: general case}
\label{subsubsec:MeasurSelec_General}

	We shall extend Theorem \ref{theo:measurable_section_t} and Proposition \ref{prop:universal_section_t} to a general context, where $\R^+$ is replaced by a more general topological space.
	The generalization of the auxiliary space $\R^+$ may be an important issue for the applications. For example, in a stochastic control problem that we shall study later,
    the auxiliary space is chosen as the space of probability measures $\Pc(E)$ on a Polish space $E$, which is also a Polish space equipped with the weak topology.\\[-8mm]

\paragraph{\small Isomorphism between $[0,1]$ and Borel space}

	The main idea of the extension is to show that some abstract space is ``equivalent'' to $\R^+$ (or $[0,1]$) in the measurable sense.

	\begin{Definition}
		\rmi A topological space is said to be a Borel space, if it is topologically homeomorphic to a Borel subset of a Polish space.\\
		\rmii Let $E$ and $F$ be two Borel spaces, $E$ and $F$ is said to be isomorphic, if there is a bijection $\varphi$ between $(E, \Bc(E))$ and $(F, \Bc(F))$ such that $\varphi$ and $\varphi^{-1}$ are both measurable.
	\end{Definition}

	Let $\varphi$ be an isomorphic bijection between Borel spaces $E$ and $F$, $\mu$ a positive measure on $(E, \Bc(E))$ and $N$ a $\mu$-null set.
    Then clearly, $\varphi(N)$ is also a null set under the imagine measure of $\mu$, it follows that for every $A \in \Bc^U(E)$, we have $\varphi(A) \in \Bc^U(F)$.
	Further, it is clear that every Polish space is  Borel space.
    But the more important result for us is the links with $[0,1]$.
    The following result is a classical one, whose proof may be found in Chapter 7 of Bertsekas and Shreve \cite{Bertsekas_1978}.
	
	\begin{Lemma} \label{lemm:Borel_Isomorph}
			Every Borel space is isomorphic to a Borel subset of the unit interval $[0,1]$. In particular, if the Borel space is uncountable, it is isomorphic to $[0,1]$.
	\end{Lemma}

	The above Lemma says that from a measure theoretic point of view, a Borel space is identical to a Borel subset of $[0,1]$.
    Then Theorem \ref{theo:measurable_section_t} holds true with slight modifications when $(\R^+, \Bc(\R^+))$ is replaced by a Borel space.
    However, since the above isomorphism is not constructive, the general measurable selection theorem given below is only an existence result. \\[-8mm]

\paragraph{\small Measurable selection theorem}

	Finally, in resume, let us give a universally measurable selection theorem in a general context, where the notation $+\infty$ in the $\R^+$ case is replaced by a cemetery point $\partial$.

	\begin{Theorem}\label{theo:measurable_section}
		Let $(\O, \Fc)$ be a measurable space, $E$ be a Borel space with $\Ec := \Bc(E)$, and $A \in \Fc \otimes \Ec$ be a measurable subset in $\O \x E$.
	 Then there exists a $\Fc^U$-universally measurable mapping $Z$ from $(\O, \Fc^U)$ into $ \big( E \cup \{\partial \}, \Kc \bigvee \big\{ \emptyset, \{\partial\} \big\} \big)$ such that
		\be
			\big\{ (\om, Z(\om)) \in \O \x E \big\} \subseteq A, ~~~\text{and}~~\big\{ \om ~: Z(\om) \in E \big \} ~=~ \pi_{\O}(A).
		\ee
	\end{Theorem}

	In the same spirit, we can also study a similar optimization problem as in Corollary \ref{Coro:Measurability_Sup_K}.
    Let $(\O, \Fc)$ and $(E, \Ec)$ be given as in Theorem \ref{theo:measurable_section}, $A \in \Fc \otimes \Ec$ and $f$ a $\Fc \otimes \Ec$-measurable function, denote, for every $\eps > 0$,
	\b*
		g(\om) := \sup \big \{ f(\o, x) ~: (\om, x) \in A \big \} &\mbox{and}& g^\eps(\om) := \big( g (\om) - \eps \big) \1_{g(\om) < \infty} + \frac{1}{\eps} \1_{g(\om) = \infty}.
	\e*

	\begin{Proposition} \label{Prop:Measurability_Sup}
		\rmi The function $g$ is $\Fc^U$-universally measurable, taking value in $\R \cup\{ \infty, -\infty \}$. Moreover, there is a $E$-valued,
        $\Fc^U$-measurable variable $Z_{\eps}$ for every $\eps > 0$ such that $\forall \om \in \pi_{\O}(A)$,
		\be \label{eq:meas_selec_sup}
			Z_{\eps}(\om) \in A^{\om}
			&\mbox{and}&
			f(\om, Z_{\eps}(\om)) ~\ge~ g^\eps(\om).
		\ee
		\rmii If $(\O, \Fc)$ is equipped with a probability measure $\P$, we can exchange ``supremum'' and ``expectation'' operators, i.e.
		\b*
			 \E^{\P} \big[ g(\om) 1_{\pi_{\O}(A)} \big]
			~=~
			 \sup \Big \{ \E^{\P} \big[ f(\om, Z(\om)) 1_{\pi_{\O}(A)} \big] ~: Z\in \Lc(\Fc) ~\mbox{s.t.}~ Z(\om) \in A^{\om}, ~\P\mbox{-a.s.} \Big \}.
		\e*
	\end{Proposition}

	\proof \rmi The measurability of $g$ is a direct consequence of Proposition \ref{prop:universal_section_t} together with the fact that $E$ is isomorphic to a Borel subset of $[0,1]$ by Lemma \ref{lemm:Borel_Isomorph}. Then
	it is enough to consider the $\Fc \otimes \Ec$-measurable set $A_{\eps} := \big \{ (\om, x) \in A ~: f(\om, x) \ge g^\eps(\om) \big \}$ in the product space and apply Theorem \ref{theo:measurable_section} to choose $Z_{\eps}$.\\
	 \rmii Given the probability measure $\P$ on $(\O, \Fc)$, then there is a $\Fc$-measurable r.v. $ Z_{\eps}$ such that \eqref{eq:meas_selec_sup} holds true for $\P$-a.e. $\om \in \pi_{\O}(A)$.
        It follows that the second assertion holds true. \qed

	\vspace{2mm}

	In Theorem \ref{theo:measurable_section} and Proposition \ref{Prop:Measurability_Sup}, $(\O, \Fc)$ is assumed to be a  measurable space without any topological structure imposed.
    It is an extension of a similar result of Bertsekas and Shreve \cite{Bertsekas_1978} Chapter 7,
    where $\O$ is assumed to be a Borel space and $\Fc$ is its Borel $\sigma$-field (see also Theorem \ref{theo:analytic_meas_selec} below).
    It may be important to consider an abstract space $(\O, \Fc)$ for many applications, for example, in the general theory of stochastic processes of Dellacherie \cite{Dellacherie_1972}.

\subsubsection{Analytic selection theorem}
\label{subsubsec:MeasSelec_analytic}

	Following the above measurable projection and selection theorems, the projection of a Borel set in the product space turns to be a universally measurable set.
	In the topological case, Souslin gave a more precise description, that is the projection set of a Borel set is an ``analytic'' set.
  A remarkable property is that the ``analytic'' sets are stable under projection, which makes the presentation simpler
  when we want to compose the projection or supremum for several times,
  as we shall see later in the gambling house model as well as the nonlinear operators framework in Section \ref{sec:DPP}.

	There are several equivalent ways to define an analytic set, we shall give one and refer to Bertsekas and Shreve   \cite{Bertsekas_1978} as well as more recent Bogachev book \cite{Bogachev} for more details.

	\begin{Definition} \label{def:analytic}
        \rmi Let $E$ be a Borel space, then a subset $B$ is an analytic set in $E$ if there is another Borel space $F$ and a Borel subset $A \subseteq E \x F$ such that $B = \pi_E(A)$.
        A subset $C \subseteq E$ is co-analytic if its complement $C^c$ is analytic.\\
        \rmii  A function $g: E \to \Rb = \R \cup \{\infty \}$ is upper semianalytic (u.s.a.) if $\{ x \in E ~: g(x) > c \}$ is analytic for every $c \in \R$.\\
        \rmiii Let $E$  be a Borel set and   $\Ac(E)$ denote the $\sigma$-field generated by all analytic subsets.
        A function  $f: E \to F$, where $F$ is a Borel set, is analytically measurable if $f^{-1}(C) \in \Ac(E)$ for every $C \in \Bc(F)$.
	\end{Definition}

		In a Borel space $E$, every Borel set is analytic, every analytic set is universally measurable, i.e. $\Bc(E) \subset \Ac(E) \subset \Bc^U(E)$.
        It follows that every upper semianalytic function is universally measurable. However, the complement of an analytic set may not be analytic and the class of analytic sets is not a $\sigma$-field.
		Nevertheless, projection and selection theorems may be extended in this context.

	\begin{Theorem} \label{theo:analytic_meas_selec}
        Let $E$ and $F$ be Borel spaces, $A$ be an analytic subset of $E \x F$, and $f : A \to \R$ be an upper semianalytic function. Define $g(x) := \sup_{(x,y) \in A} f(x,y)$.\\
		\rmi The projection set $\pi_E(A)$ is an analytic subset in $E$.\\
		\rmii There exists an analytically measurable function $\varphi : \pi_E(A) \to F$ such that $(x,\varphi(x)) \in A$, for every $x \in \pi_E(A).$\\
		\rmiii The function $g: \pi_E(A) \to \Rb = \R \cup \{ \infty \}$ is upper semianalytic.\\
		\rmiv For every $\eps > 0$, there is an analytically measurable function $\varphi_{\eps}: \pi_E(A) \to F$
        such that $f(x,\varphi_{\eps}(x)) \ge g^{\eps}(x) := \big( g(x) - \eps \big) 1_{g(x) < \infty} + \frac{1}{\eps} 1_{g(x) = \infty}$ for every $x \in \pi_E(A)$.
	\end{Theorem}

    The assertion \rmi follows directly from the definition of the analytic set, \rmiii and \rmiv are standard sequences of \rmi and (ii).
    The assertion \rmii is more essential, which cannot be deduced by the projection capacity approach.
    The main idea of the proof is to use the fact that an analytic set is a continuous image of the Baire space (the set of all infinite sequences of natural numbers)
    and then to explore the topological properties of the Baire space.
    Let us refer to Chapter 7 of Bertsekas and Shreve \cite{Bertsekas_1978} or Chapiter 2 of Bogachev \cite{Bogachev} (vol 2) for a complete technical proof.


\section{Discrete time dynamic programming}
\label{sec:DPP_discrete}

    It is classical to use the measurable selection theorem to deduce the dynamic programming principle (DPP) or time consistence property in discrete time models,
    see for example the gambling house models studied in Dellacherie \cite{Dellacherie_MaisonJeux}, Maitra and Sudderth \cite{MaitraSudderth},
    and also the control problem in Bertsekas and Shreve \cite{Bertsekas_1978} (for which a brief version can also be found in the expository paper of  Bertsekas and Shreve \cite{BertsekasShreve}).

\subsection{Probability kernel, composition and disintegration}

	We here recall the definition as well as some facts on probability kernel, the composition and disintegration of the probability measures, for which our main reference is Dellacherie and Meyer \cite{Dellacherie_Meyer_1988}.\\ [-8mm]

\paragraph{Probability kernel}

	The probability kernel can be viewed as a probabilistic generalization of the concept of function: instead of associating with each input value a deterministic output value, one chooses a random output value or equivalently a probability distribution.

	\begin{Definition}
		Let $(X, \Xc)$ and $(Y, \Yc)$ be two measurable spaces, a probability kernel from $(X, \Xc)$ to $(Y, \Yc)$ (or simply say from $X$ to $Y$) is a map $N : X \x \Yc \to [0, 1]$ such that \rmi for any $B\in \Yc$ the function $x \mapsto N(x,B)$ is $\Xc$-measurable; \rmii for any $x \in X$, $N(x,\cdot)$ is a probability measure on $(Y, \Yc)$.
	\end{Definition}

	Let $N$ be a probability kernel from $(X, \Xc)$ to $(Y, \Yc)$, $f$ a bounded $\Yc$-measurable function, we define a bounded function $Nf : X \to \R$ by $Nf(x) := \int_Y\,f(y) N(x,dy)$.
	It is then clear that $Nf$ is a $\Xc$-measurable function.
	Moreover, let $(f_n)_{n \ge 1}$ be an increasing sequence of bounded measurable functions on $X$ such that $f_n \uparrow f$,
	then $N(\lim f_n) = \lim N f_n$.
	Another important property of the probability kernel is that it can be extended uniquely to its univerally completion.
	
	\begin{Lemma} \label{lemm:universal_kernel}
		Let $(X, \Xc)$ and $(Y, \Yc)$ be two measurable spaces, $\Xc^U$ and $\Yc^U$ be the universal completions of $\Xc$ and $\Yc$.
		Then $N$ extends uniquely to a kernel from $(X, \Xc^U)$ to $(Y, \Yc^U)$.
	\end{Lemma}

\paragraph{Conditional probability distribution (c.p.d.)}

		An important example of probability kernel is the so-called conditional probability distribution (c.p.d.), or the disintegration of a probability measure.
		Let $\P$ be a probability measure on a measurable space $(X, \Xc)$, $\Uc$ be a sub-$\sigma$-field of $\Xc$, a conditional probability distribution of $\P$ w.r.t $\Uc$ is a family of probability measures $(\P_x)_{x \in X}$ such that $x \mapsto \P_x$ is $\Uc$-measurable and for any $C\in \Xc$, $\P_x(C)=\P(C|\Uc)(x)$ for $\P$-almost every $x \in X$.
		
		It is clear that the map $(x, A) \mapsto \P_x(A)$ for all $(x, A) \in (X, \Xc)$ defines a a probability kernel from $(X, \Uc)$ to $(X, \Xc)$.
    We further notice that when $X$ is a Polish space, $\Xc = \Bc(X)$ its Borel $\sigma$-field, $\Uc$ a sub $\sigma$-field of $\Bc(X)$ and $\P$ a probability measure defined on $(X, \Xc)$, the conditional probability distribution of $\P$ w.r.t $\Uc$ exists.

    Let $(X, \Xc)$ be a measurable space equipped with a probability measure $\P$, assume the existence of the c.p.d. $(\P_x)_{x \in X}$ of $\P$ w.r.t. $\Uc$, suppose in addition that $\Uc$ is countably generated,
    then there is a $\P$-negligible set $N$ such that for any $x\in N^c$ and $U\in \Uc$, $x \in U \Longrightarrow \P_x(U)=1$.
    In other words, for $x\in N^c$, $\P_x$ and the Dirac measure $\delta_{x}$ are the same restricted on $\Uc$.
    Further, by setting $\P_x := \delta_{x}, ~\forall x \in N$, we obtain a particular family of c.p.d.,
    which is called a regular conditional probability distribution (r.c.p.d.) following the terminology of Stroock and Varadhan \cite{Stroock_1979}.\\[-8mm]

\paragraph{Composition of probability kernels}
	Given a probability measure as well as a probability kernel, or two probability kernels, we can define their compositions (or concatenations).
	Let $\mu$ be a probability measure on $(X, \Xc)$ and $N$ a probability kernel from $X$ to $Y$, the composition $\mu \ox N : \Yc \to [0,1]$ is defined by $\mu \ox N (A) := \int_X \mu(dx) N(x, A)$.
	One can easily verify that $\mu \ox N$ is a probability measure on $(Y, \Yc)$.
	Moreover, $\langle \mu, N f \rangle = \langle \mu \ox N, f \rangle$ for every bounded measurable function defined on $Y$.
	Let $M$ be a probability kernel from $(X,\Xc)$ to $(Y,\Yc)$, $N$ a kernel from $(Y,\Yc)$ to $(Z,\Zc)$,
	we can then define a probability kernel $M \ox N$ from $X$ to $Z$ by $M \ox N (x,A) := \int_Y M(x, dy) N(y,A)$.

\subsection{Gambling house}
	A typical example, widely studied in the 70'st by Dellacherie \cite{Dellacherie_MaisonJeux}, Maitra and Sudderth \cite{MaitraSudderth} and many others,
    is the gambling house model which involves a Borel space $E$ and  the Borel space $F=\Pc(E)$ of probability measures on $E$ equipped with the weak convergence topology.

	A {\em gambling house} is an analytic subset $J$ of $E \times \Pc(E)$ such that the section set $J_x := \{\mu ~: (x,\mu) \in J \}$ is nonempty for every $x \in E$;
	this name is motivated by the interpretation of the section set $J_x$ as the collection of distributions of gains available to a gambler having a wealth $x$.
	The maximal expected gain $Jf : E \to \R^+$ associated with an upper semianalytic reward function $f: E \to \R^+$ is defined  by
    \be \label{maxgain}
   		Jf(x) ~:=~ \sup_{\mu\in J_x} \langle \mu, f \rangle.
    \ee
    It follows by Theorem \ref{theo:analytic_meas_selec} that $Jf$ is also an upper semianalytic function.
	We notice that when the sections $J_x$ is compact, then $f\mapsto J f(x)$ is in fact a (functional) capacity, in spirit of Proposition \ref{supmeasure} and Remark \ref{rem:capacity_topo}.

	We can also define the compositions of gambling houses.
	Recall that $\Ac(E)$ denotes the $\sigma$-field generated by all analytic sets in $E$, which is included by the universal $\sigma$-field $\Bc^U(E)$.
	Let $J$ be a gambling house, a probability kernel $\lambda$ from $(E, \Bc(E))$ to $(E, \Ac(E))$, is said to be $J$-admissible if $(x,\lambda(x,\cdot)) \in J$ for every $x \in E$.
	Given a probability $\mu \in \Pc(E)$, and two gambling houses $J$ and $K$, we define $\mu \ox K \subseteq \Pc(E)$ by
	\b*
		\mu \ox K
		&:=&
		\Big\{ \mu \ox \lambda ~: \lambda ~ K\mbox{-admissible kernel} \Big\},
	\e*
	and a composed gambling house $J \ox K$ by
	\b*
		J \ox K &:=& \big\{ (x, \lambda) ~: \lambda \in \mu \ox K, ~  (x,\mu) \in J \big\}.
	\e*
	Clearly, by Lemma \ref{lemm:universal_kernel}, the composition $\mu \ox \lambda$ is uniquely defined given a $K$-admissible (analytic) kernel $\lambda$.
	Similarly, for $n$ gambling houses $J^1, \cdots, J^n$, their composition gambling house is given
	by $J^1 \ox \cdots \ox J^n := J^1 \ox ( \cdots \ox (J^{n-1} \ox J^n))$.
	
	\begin{Proposition} \label{prop:maison_jeux}
		Let $\Phi : E \to \R^+$ be an upper semianalytic reward function, $\mu \in \Pc(E)$, $J$ and $(J^k)_{1 \le k \le n}$ be a sequence of gambling houses. Then\\
		\rmi $J(\Phi)$ is stille upper semianalytic and 
            \be \label{eq:mJ}
                m \ox J(\Phi) := \sup_{\mu \in m \ox J} \langle \mu, \Phi \rangle 
                ~~=~~ 
                \int_E J \Phi(x) m(dx).
            \ee
		\rmii We have $J^1 \ox \cdots \ox J^n (\Phi) = J^1( \cdots J^{n-1}( J^n(\Phi)))$.
	\end{Proposition}
    \proof 
        Since $\Phi$ is upper semianalytic, then the map $(x, \mu)\in E \x \Pc(E) \mapsto \langle \mu, \Phi \rangle \in \R$ is also upper semianalytic.
        It follows by \eqref{maxgain} together with the measurable selection theorem (Theorem \ref{theo:analytic_meas_selec}) that $J(\Phi)$ is upper semianalytic and for every $\eps > 0$, there is $J$-admissible analytic kernel $\lambda_{\eps}$ such that $\lambda_{\eps}(x, \Phi) \ge J \Phi(x) - \eps$ for all $x \in E$.
        Therefore, 
        \b*
            m \ox J(\Phi) 
            &:=& 
            \sup_{\mu \in m \ox J} \langle \mu, \Phi \rangle
            ~\ge~
            m \ox \lambda_{\eps} (\Phi)\\
            &=&
            \int_E \big( \lambda_{\eps}(x, \Phi) \big) m(dx)
            ~\ge~
            \int_E J \Phi(x) m(dx) - \eps.
        \e*
        Further, for every $\mu \in m \ox J$, there is $J$-admissible kernel $\lambda$ such that $\mu = m \ox \lambda$.
        It follows that
        \b*
            \langle \mu, \Phi \rangle 
            &=&
            \int_E \lambda(x, \Phi) m(dx)
            ~\le~
            \int_E J \Phi(x) m(dx).
        \e*
        We then conclude the proof for \rmi by the arbitrariness of $\eps > 0$ and $\mu \in m \ox J$.
        Finally, by exactly the same arguments, we can easily prove \rmii.
    \qed

	\begin{Remark}
	{\rm 
        The result in Proposition \ref{prop:maison_jeux} is in fact the dynamic programming principle when the gambling house is interpreted as a discrete time control problem.
        Generally, the equality in \eqref{eq:mJ} is estabilished by two inequalities using different arguments.
        For the first inequality, the essential is to choose a family of $\eps$-optimal $\mu$ in \eqref{maxgain} in a measurable way, so that their composition with $m$ lies in the set of arguments of the previous optimization problem in \eqref{eq:mJ}.
        For the second inequality, the essential is to be able to decompose any $\mu \in m \ox J$ in the way $\mu = m \ox \lambda$ such that $(x, \lambda(x,\cdot)) \in J$ for $m$-almost surely $x \in E$.
    }
	\end{Remark}

\section{Continuous time dynamic programming}
\label{sec:DPP}

  We shall now extend the discrete time dynamic programming to a continuous time context.
  We first introduce a family of operators on the functional space of the canonical space of c\`adl\`ag trajectories as well as its extension spaces, which are indexed by stopping times.
	Using measurable selection theorems, we show how the family of operators admits a time consistence property under appropriate conditions.
	In particular, when the operators are defined by control problems, the time consistence turns to be the dynamic programming principle of the control problem.

\subsection{The canonical space of c\`adl\`ag trajectories}
\label{subsec:canonical_space}

	Let $E$ be a Polish space, we shall consider the canonical space $\Om := D(\R^+, E)$ of all $E$-valued c\`adl\`ag paths on $\R^+$.
    Some basic properties is given in the following, where our main reference is Chapter 3 of Ethier and Kurtz \cite{Ethier_Kurtz} as well as Chapter IV of Dellacherie and Meyer \cite{Dellacherie_Meyer_1978}.
    First, $\Om$ equipped with the Skorokhod topology is a Polish space, and its Borel $\sigma$-field is also the canonical $\sigma$-field generated by the canonical process $X$ defined by $X_t(\om) := \om_t, ~\forall \om \in \Om$,
    i.e. $\Bc(\Om) = \Fc_{\infty} ~:=~ \sigma(X_t~: t \ge 0).$\\[-8mm]

\paragraph{\small Operations on the trajectories}
    There is a family of operators from canonical space $\Om$ to itself, allowing  different transformations on  paths, such as stopping and concatenating (pasting).
    Let us first define two sets of paths, given $(\w, t) \in \Om \x \R^+$,

- $\Dc_{\w}^t := \{\om ~: X_t(\om)=X_t(\w)\}$ is the set of the paths which coincide at time $t$ with $\w$,
	
- $\Dc_{(\w,t)} := \{\om ~: X_s(\om)=X_s(\w), ~\forall s \le t\}$ is the set of the paths which coincide with $\w$ on interval $[0,t]$.

    Let $t \in \R^+$, we then define three operators on the paths:\\
\rmi Stopping  operator $a_t : \Om \to \Om$:\\
\centerline{$ a_t(\om) := [\om]_t := \om_{t \land \cdot}\Longleftrightarrow X_s (a_t(\om)) = X_{s \land t}(\om), ~ \forall s \ge 0.$}
The image of $a_t$ is the collection of all paths stopped at time $t$; on the other hand, it is clear that $(a_t)^{-1}([\w]_t) = \Dc_{(\w,t)}$, $\forall \w \in \Om$.\\
\rmii Predictable concatenation  operator $\ox_{t^-} : \Om \times \Om \to \Om$: \\
\centerline{$(\w \ox_{t^-} \om)(s) := \w_s, \> \mbox{on}\>\{s < t\} ,\quad \mbox{and} \quad (\w \ox_{t^-} \om)(s) :=\om_s \quad \mbox{on} \quad \{s \ge t\}$.~~~~}\\
\rmiii Optional concatenation  operator $\ox_{t} : \{ (\w, \om) ~: \om \in \Dc_{\w}^t\} \to \Om$: \\
\centerline{$(\w \ox_{t} \om)(s) := \w_s, \> \mbox{on}\>\{s \leq t\} ,\quad \mbox{and} \quad (\w \ox_{t} \om)(s) :=\om_s \quad \mbox{on} \quad \{s >t \}$.~~~~}\\
We notice that the predictable concatenation operator $\ox_{t^-}$ loses the information given by $X_t(\w)=\w_t$, this is the main reason we consider the optional concatenation $\ox_{t}$.
In the definition of $\ox_{t}$, the condition $\om \in \Dc_{\w}^t$ ensures that $\w \ox_{t} \om$ is right continuous and hence lies in $\Om$.
In practice, we usually impose a probability measure on $\Om$ to ensure that $\om \in \Dc_{(\w,t)} \subset \Dc_{\w}^t$ almost surely, and hence the two concatenation operators become ``almost'' the same. \\[-8mm]

%
\paragraph{\small Canonical filtration}
As usual, we define the canonical filtration $\F = (\Fc_t)_{t \ge 0}$, generated by the canonical process, by $\Fc_t := \sigma(X_s, s \le t)$.
In particular, $\Fc_{\infty} = \bigvee_{t \ge 0} \Fc_t$ is the Borel $\sigma$-field of the Polish space $\Om$.
We notice that $\F$ is not right continuous; a right continuous filtration $\F^+ = (\Fc^+_t)_{t \ge 0}$ can be defined by $\Fc^+_t := \bigcap_{\eps > 0} \Fc_{t+\eps}$.\\
\rmi For every $t \in \R^+$, the $\sigma$-field $\Fc_t$ is generated by the stopping operator $a_t$ in the sense that $\Fc_t = a_t^{-1}(\Fc_{\infty})$.
	In particular any  random variable $Y$ is $\Fc_t$-measurable if and only if  $Y(\om)=Y(a_t(\om)),~ \forall \om \in \Om$.\\
\rmii The map $(t,\om) \in \Om \x \R^+ \mapsto a_t(\om) \in \Om$ is progressively measurable w.r.t. $\F$.
	Moreover, let $Z$ be a $\Bc(\R^+) \ox \Fc_{\infty}$-measurable process, then $Z$ is $\F$-progressively measurable if and only if $Z_t(\om)=Z_t(a_t(\om))$ for every $(t,\om) \in \R^+ \x \Om$.\\
\rmiii A random variable $\tau : \Om \to \R^+ \cup \{+\infty \}$ is a $\F$-stopping time if and only if for all $\om,~ \om' \in \Om$,\\
	\centerline{$\{ \tau(\om) \le t ~: a_t(\om) = a_t(\om') \},~ \forall t \ge 0 ~~\Rightarrow~~ \{ \tau(\om) = \tau(\om')  \}.$}
Let us denote by $\Tc$ the collection of all $\F$-stopping times taking value in $[0, \infty)$. \\
\rmiv Given a finite $\F$-stopping time $\tau \in \Tc$, a random variable $Y$ is $\Fc_{\tau}$-measurable if and only if there exists a measurable process $Z$ such that $Y(\om) = Z \big(\tau(\om), a_{\tau(\om)}(\om) \big)$.
	It follows that $\Fc_{\tau} = \phi^{-1}( \Bc(\R^+) \ox \Fc_{\infty})$ with $\phi(\om) := (\tau(\om), a_{\tau(\om)}(\om)) \in \R^+ \x \Om,~ \forall \om \in \Om$;
    and in particular, $\Fc_{\tau}$ is countably generated since $\Bc(\R^+) \ox \Fc_{\infty}$ is countably generated. \\[1mm]
\rmv {\em Test de Galmarino}  \cite{Dellacherie_Meyer_1978}(Th 103): Let $S$ and $T$ be two $\F$-stopping times such that $S\leq T$.
There exists a function $U(\w,\om)$ defined $\Om \times \Om$ taking values in $\R^+\cup \{\infty \}$ and $\Fc_S \times \Fc_{\infty}$-measurable,  such that

$-$ $U(\w,\om)=+\infty$ when $S(\w)=+\infty$, or when $S(\w)<+\infty$ and $\om \not \in \Dc_{(\w,S(\w))}$.

$-$ $U(\w,\om)=T(\w\ox_{S(\w)} \om)$ when $S(\w)<+\infty$ and $\om \in \Dc_{(\w,S(\w))}$.

$-$ For any $\w$, $U(\w,.)$ is a $\F-$stopping time greater than $S(\w)$. This follows in fact immediately from the characterization \rmiii of the stopping times.\\
\rmvi The predictable $\sigma$-field $\Pc$ is the $\sigma$-field defined on $\Om \times \R^+$, generated by all progressively measurable, left continuous processes; while the optional $\sigma$-field $\Oc$ is generated by all right continuous, left limited progressive processes. On the canonical space, we then have

$-$ The optional $\sigma$-field $\Oc$ is generated by $\Pc$ and $[X_t(\om)]$, where $[x]$ denotes the constant path equal to $x$.

$-$ The map $(\w,t, \om)\mapsto \w \ox_{t^-} \om$ is measurable from $\Pc\otimes \Fc_{\infty}$ into $\Fc_{\infty}$.

$-$ The graph of $(\Dc_{(\w,t)})_{(\w,t) \in \Om \x \R^+}$ defined as $[[\Dc]] :=\{(\w,t, \om) ~: \om \in \Dc_{(\w,t)}\}$ lies in $\Oc \otimes \Fc_{\infty}$.

$-$ The map $(\w,t,\om)\mapsto \w \ox_{t} \om$ restricted to $[[\Dc]]$  is measurable from $\Oc \otimes \Fc_{\infty}$ into $\Fc_{\infty}$.\\[1mm]
To conclude, we notice that an important class of stopping times is the hitting times in canonical space.
These hitting times can be considered as the debut of sets, whose universally measurability property is also a key tool in Proposition \ref{prop:measurability_debut} to deduce the measurable section result.
Following Dellacherie \cite{DellacherieAnalytique}, we give some more precise measurability properties of the hitting times.
\begin{Corollary} [Measurability of hitting times]\label{analytic debut}
    Let  $A$ be an analytic subset of $E$.
	The debut and the hitting time of $A$ are defined respectively by
	\b*
        D_A(\o) ~:=~ \inf \{ t\geq 0 ~: X_t(\o)\in A\}
        ~~ \mbox{and}~~
        T_A(\o) ~:=~ \inf \{ t > 0 ~: X_t(\o)\in A\}.
    \e*
\rmi First, $D_A$ and $T_A$ are coanalytic random variables. \\
\rmii Suppose that $A$ is open, then $D_A$ and $T_A$ are Borel measurable random variables and stopping times w.r.t. the right continuous filtration $\F^+$. \\
\rmiii Suppose that $A$ is closed, then the same properties hold true for the variables $Z_A$ and $S_A$,
    where $Z_A(\o) := \inf \{ t\geq 0 ~: X_t(\o) ~\mbox{or}~ X_{t_-}(\o) \in A\}$ and $S_A(\o) := \inf \{ t > 0 ~: X_t(\o)\> \text{or}\>X_{t_-}(\o)\in A\}$.
    More precisely, $Z_A$ is in fact a $\F$-stopping time.\\[-8mm]
\end{Corollary}
\paragraph{\small Conditioning and concatenation on the canonical space}
Let $\Pc(\Om)$ denote the space of all probability measures on $\Om$, which is also a Polish space under the weak convergence topology.
We can then define the conditioning (or disintegration) and concatenation of probability measures in $\Pc(\Om)$ as discussed in Dellacherie, Meyer \cite{Dellacherie_Meyer_1978}, or in Stroock and Varadhan \cite{Stroock_1979}.
In the canonical space case, the operators defined on the paths allow us to give an intuitive description of the regular conditional operator, very similar to Markovian kernels.
A short resume of the weak convergence topology as well as conditioning and concatenation of probability measures is provided in Appendix.
\\
\rmi Let $\P \in \Pc(\Om)$ and $\tau \in \Tc$ be a $\F$-stopping time taking value in $[0,\infty)$.
Since $\Fc_{\tau}$ is countably generated, there exists a family of regular conditional probability distribution (r.c.p.d.) $(\P^{\tau}_{\w})_{\w \in \Om}$ of $\P$ with respect to $\Fc_{\tau}$.
In particular, for every $\w \in \Om$, we have $\P^{\tau}_{\w}(\Dc_{(\w,\tau(\w))})=1$ and $\P^{\tau}_{\w}(A)= \1_A(\w), \forall A\in \Fc_{\tau}$.\\
%
%
%
%
%
%
%
\rmii On the other hand, suppose that $\P$ is a probability measure on $(\O, \Fc_{\tau})$ and $(\Q_{\w})_{\w\in \O}$ is a probability kernel from $(\Om,\Fc_{\tau})$ to $(\Om,\Fc_{\infty})$
    such that $\Q_{\w} (\Dc_{\w}^{\tau(\w)}) = 1, ~\forall \w \in \Om$.
    By composition of probability kernels, there is a unique concatenated probability measure $\P \ox_{\tau} \Q_{\cdot}$ on $(\Om, \Fc_{\infty})$, defined by
	\be \label{eq:concatenation_tau}
		\P \ox_{\tau}  \Q_{\cdot}(A):=\int \P(d\w)\int \1_A(\w \ox_{\tau(\w)}\om)\Q_{\w}(d\om).
	\ee
    In particular, we have $\P \ox_{\tau} \Q_{\cdot}(A) = \P(A),~ \forall A \in \Fc_{\tau}$ and $(\delta_{\w} \ox_{\tau(\w)} \Q_{\w})_{\w\in \O}$ is a family of r.c.p.d. of $\P \ox_{\tau}  \Q_{\cdot}$ w.r.t $\Fc_{\tau}$.
\begin{Remark}
    {\rm  \rmi For canonical space of real valued continuous paths $C([0,T], \R^d)$, a detailed presentation of the conditionning and concatenation is given in Stroock and Varadhan \cite{Stroock_1979}.\\
    \rmii It is shown in Blackwell and Dubins \cite{BlackDubins} that when $\Om$ is the canonical space of real valued continuous paths and $\F$ the canonical filtration,
    there exists no regular conditional distribution w.r.t. $\Fc_t^+=\bigcap_{\eps>0} \Fc_{t+\eps}$.
	A detailled exposition of questions related to this subject may be found in  Bogachev's Measure Theory book (Vol 2, Chap. IX, X)\cite{Bogachev}.
	Another detailed discussion with many applications to statistics may be found in Chang and Pollard \cite{ChangPollard}.}
\end{Remark}

\subsection{A family of time consistent nonlinear operators}
\label{subsec:TC}

\subsubsection{Markov property for the historical process }

    Before introducing the notion de time consistent nonlinear operators, we briefly present the case of linear operators, generally known as (non-homogeneous) semi-group.
    The time consistency is nothing else than the strong Markov property.
    Let us introduce the so-called historical process $(Y_t)_{t \ge 0}$ defined, from $\Om$ to $\Om$, by\\[2mm]
    \centerline{  $Y_t(\om)=a_t(\om), \quad$ for any $t\in \R^+$ and any $\omega \in \O $. }\\[2mm]
 It is clear that the filtration generated by the processes $Y$ is the same as that generated by the canonical process $X$.
 In particular, $\Fc_t = \Fc^Y_t=\sigma(Y_s; s\leq t)= \sigma(Y_t)$, and $\Fc_{\infty} = \Fc^Y_\infty=\sigma(Y_s; s<\infty )= \sigma(Y^-_\infty)$.
 Further, the progressively measurable adapted processes are given by $\xi(t,Y_t(\o))$ where $\xi(t,\w)$ is a measurable function on $\Omega \times R^+$.

    Suppose that $(\Q_{\w,t})_{(\w,t) \in \Om \x \R^+}$ is a given probability kernel from $\Om \x \R^+$ on $(\Om, \Fc_{\infty})$, satisfying
    \rmi  The initial condition: $\Q_{\w,t}(Y_t = a_t(\w))= 1, ~\forall (\w,t) \in \Om \x \R^+$.
    \rmii Time consistency: for every $(\w_0, t_0) \in \Om \x \R^+$ and every finite $\F-$stopping time $\tau\geq t_0$,
    $(\Q_{\om,\tau(\om)})_{\om \in \Om}$ is a r.c.p.d. of $\Q_{\w_0,t_0}$ w.r.t. $\Fc_{\tau}$.
    In particular, we see that $ \Q_{\w,t} = \Q_{[\w]_t,t} = \Q_{Y(\w),t}$.
    Then $Y$ is a strong Markov process with Markovian semigroup $(\Q_{\w,t})_{(\w,t) \in \Om \x \R^+}$.
	More generally, let $\tau$ and $\sigma$ be two finite stopping times such that $\tau \leq \sigma$.
    Then, for any positive $\Fc_{\infty}$-measurable function $\xi$, we also have
    \be \label{DDPMarkov}
        \E^{\Q_{(Y_\tau(\om),\tau(\om))}} \big[ \xi(Y^-_{\infty}) \big]
        &=&
        \E^{\Q_{(Y_\tau(\om),\tau(\om))}} \big[ \E^{\Q_{(Y_\sigma,\sigma)}}[\xi(Y^-_{\infty})] \big].
    \ee
	That is the classical version of the dynamic programming principle etablished in the sequel.
	\begin{Remark}
    {\em An example of such families of probability measures is that induced by a family of diffusion processes $(X^{\w,t})_{(\w,t) \in \Om \x \R^+}$, defined by $X^{\w,t}_{\theta} := \w_{\theta}$ when $\theta \le t$, and
    \b*
    		X^{\w,t}_{\theta} := \int_t^{\theta} \mu(s, [X^{\w,t}]_s) ds + \int_t^{\theta} \sigma(s, [X^{\w,t}]_s) dW_s,
    		~~\mbox{when}~ \theta > t,
    \e*
    where $E = \R^d$ and $\mu: \R^+ \x \Om \to \R^d$, $\sigma : \R^+ \x \Om \to S_d$ is the diffusion coefficient.}
	\end{Remark}

%

	
\subsubsection{A family of time consistent nonlinear operators}
\label{subsubsec:TC}

	Let us denote by $\Ac_{usa}(\Om)$ the collection of all upper semianalytic (u.s.a.) functions bounded from below defined on the Polish space $\Om$.
	Given a $\F-$stopping time $\tau$, we denote by $\Fc_{\tau}^U$ the universally completed $\sigma$-field of $\Fc_{\tau}$
	and by $\Ac_{\tau}^U(\Om)$  the collection of all $\Fc_{\tau}^U$-measurable functions in $\Ac_{usa}(\Om)$.
	We shall consider a family of nonlinear operators associated with a class of probability measure families $(\Pc_{t,\om})_{(t,\om) \in \R^+ \x \O}$:
	\be \label{eq:defEc}
		\Ec_{\tau} \big[\xi \big](\om)
		&:=&
		\sup \big\{ \E^{\P} \big[ \xi \big]\> : \P \in \Pc_{\tau(\om),\om} \big\},~~\forall ~\mbox{finite}~ \F \mbox{-stopping time}~ \tau.
	\ee

	The family $(\Pc_{t,\om})_{(t,\om) \in \R^+ \x \O}$ can be considered as a family of section sets of a subset in $\R^+ \x \Om \x \Pc(\Om)$.
	Equivalently, we consider its graph
	\be \label{eq:Pc_graph}
		[[\Pc]] :=
		\big \{ (t, \om, \P) ~: (t,\om) \in \R^+ \x \Om,~ \P \in \Pc_{t, \om} \big\} .
	\ee
	Suppose that $[[\Pc]]$ is an analytic set in the Polish space $\R^+ \x \Om \x \Pc(\Om)$.
	Moreover, we assume the progressive measurability, i.e. for every $(t,\w) \in \R^+ \x \Om$,  $\Pc_{t,\w} $ is not empty and $\F$-adapted with support $\Dc_{(\w,t)}=\{\om ~: a_t (\om)= a_t(\w)\}$.
	In other words, for any $(t,\w) \in \R^+ \x \Om$,
	\be \label{eq:Prop_Pc}
		\Pc_{t,\w} = \Pc_{t, [\w]_t}\quad~\mbox{and}~ \quad  \P\big( \Dc_{(\w,t)}\big) = 1, ~~\forall \> \P \in \Pc_{t,\w}.
	\ee

	\begin{Lemma} \label{lemm:measurability_TC}
		Let $(\Pc_{t,\om})_{(t,\om) \in \R^+ \x \O}$ be given above, $\tau \in \Tc$ and $\xi \in \Ac_{usa}(\Om)$.
		Then $\Ec_{\tau}(\xi) \in \Ac^U_{\tau}(\Om)$.
		In particular, $\Ec_{\tau}$ is an operator from $\Ac_{usa}(\O)$ to $ \Ac_{\tau}^U(\Om) \subset \Ac_{usa}(\O)$.
	\end{Lemma}
	\proof For every positive upper semianalytic (u.s.a.) function $\xi$, the map $\P \mapsto \E^{\P}[\xi]$ is also u.s.a by Corollary 7.48.1 of \cite{Bertsekas_1978}. We then introduce
	\b*
		V(t, \om) &:=& \sup_{(t,\om, \P) \in [[\Pc]]}  \E^{\P}[\xi],~~~ \forall (t,\om) \in \R^+ \x \Om,
	\e*
	which is also u.s.a. on $\R^+ \x \Om$ from Theorem \ref{theo:analytic_meas_selec}.
	We then conclude the proof by the fact that $\Ec_{\tau}(\xi)(\om) = V(\tau(\om), \om)$.\qed

	\vspace{2mm}

	Our main objective is to derive a time consistency property for operators $(\Ec_{\tau})_{\tau \in \Tc}$ indexed by stopping times as a reformulation of the dynamic programming principle.
	From another point of view, it is formally a permutation property between the supremum and the expectation as in Proposition \ref{Prop:Measurability_Sup}.
	We notice that a closed framework is proposed in Nutz and van Handel \cite{NutzHandel} on canonical space of real valued continuous paths.
	The following assumptions generalize the notion of Markov kernel for the historical process $Y$.

	\begin{Assumption} \label{assum:stability}
		Let $(t_0,\om_0) \in \R^+ \x \Om$ be arbitrary, $\tau$ be an arbitrary stopping time taking value in $[t_0, \infty)$ and $\P \in \Pc_{t_0,\om_0}$.\\
		\rmi({\bf Stability by conditioning})
		There is a family of regular conditional probability measures $(\P_{\om})_{\om \in \Om}$ of $\P$ w.r.t. $\Fc_{\tau}$ such that $\P_{\om} \in \Pc_{\tau(\om), \om}$ for $\P$-almost every $\om \in \Om$.\\
		\rmii({\bf Stability by concatenation})
		Let $(\Q_{\om})_{\om \in \Om}$ be such that $\om \mapsto \Q_{\om}$ is $\Fc_{\tau}$-measurable and $\Q_{\om} \in \Pc_{\tau(\om), \om}$ for $\P$-a.e. $\om \in \Om$, then $\P \ox_{\tau} \Q_{\cdot} \in \Pc_{t_0, \om_0}$.
	\end{Assumption}

	\begin{Theorem} \label{theo:TC}
		Suppose that $[[\Pc]]$ is analytic in $\R^+ \x \Om \x \Pc(\Om)$, the condition \eqref{eq:Prop_Pc} and Assumption \ref{assum:stability} hold true.
        Then for every stopping times $\tau \le \sigma \in \Tc$, we have the following time consistence property:
		\be \label{eq:TC}
			  \Ec_{\tau} [\xi] &=& \Ec_{\tau} \big[ \Ec_{\sigma} [\xi] \big], ~~~ \forall \xi \in \Ac_{usa}(\Om).
		\ee
	\end{Theorem}

	\proof Let $\xi \in \Ac_{usa}(\O)$, $\om_0 \in \Om$ and $\P \in \Pc_{\tau(\om_0),\om_0}$, then following the stability assumption (Assumption \ref{assum:stability}),
	there is a family of conditional probability measures $(\P^{\sigma}_{\om})_{\om \in \O}$ of $\P$ w.r.t. $\Fc_{\sigma}$ such that $\P^{\sigma}_{\om} \in \Pc_{\sigma(\om), \om}$ for $\P$-a.e. $\om \in \Om$.
	It follows that
	\be \label{eq:TC_le}
		\E^{\P} ~ \big[ \xi \big]
		&=&
		\E^{\P} ~\Big[ \E^{\P^{\sigma}_{\om}} \big[  \xi \big] \Big]
		~~\le~~
        \E^{\P} ~ \Big[  \Ec_{\sigma} [\xi]  \Big].
	\ee

	Next, for every $\P \in \Pc_{\tau(\om_0),\om_0}$ and $\eps > 0$, denote $\Ec^{\eps}_{\sigma}[\xi] := (\Ec_{\sigma}[\xi] - \eps )1_{\Ec_{\sigma}[\xi] < \infty} + \frac{1}{\eps} 1_{\Ec_{\sigma}[\xi] = \infty}$.
	It follows by Proposition \ref{Prop:Measurability_Sup} that we can choose a family of probability $(\Q^{\eps}_{\om})_{\om \in \O}$ such that $\om \mapsto \Q^{\eps}_{\om}$ is $\Fc_{\sigma}$-measurable and
	\b*
		\Q^{\eps}_{\om} \in \Pc_{\sigma(\om), \om}, ~~~ \E^{\Q^{\eps}_{\om}}[\xi] \ge \Ec^{\eps}_{\sigma}[\xi](\om) &\mbox{for}& \P\mbox{-a.e.}~\om \in \Om.
	\e*
	Then $\P \ox_{\sigma} \Q^{\eps}_{\cdot} \in \Pc_{\tau(\om_0), \om_0}$ by the stability assumption (Assumption \ref{assum:stability}), which implies that
	\be \label{eq:TC_ge}
		\Ec_{\tau} \big[ \xi \big] ~~\ge~~ \E^{\P \ox_{\sigma} \Q^{\eps}_{\cdot}} \big[  \xi \big] ~\ge~ \E^{\P}~ \big[ \Ec^{\eps}_{\sigma} \big[\xi \big] \big].
	\ee
	We then conclude the proof by the arbitrariness of $\P \in \Pc_{\tau(\om_0),\om_0}$ and $\eps >0$ in \eqref{eq:TC_le} as well as in \eqref{eq:TC_ge}.
    \qed

 \vspace{2mm}

 Similar to Corollary 2.5 of Neufeld and Nutz \cite{NeufeldNutz}, it is easy to see that the condition in Assumption \ref{assum:stability} is stable under intersection.
 \begin{Proposition}
 		Suppose that there are two families $\Pc^1$ and $\Pc^2$ satisfying Assumption \ref{assum:stability}, $[[ \Pc^1]]$ and $[[ \Pc^2 ]]$ are both analytic.
 		Then $\Pc^0 := \Pc^1 \cap \Pc^2$ satisfies also Assumption \ref{assum:stability} and $[[\Pc^0]]$ is also analytic.    \\[-6mm]
 \end{Proposition}

\paragraph{Time consistency of dynamic risk measures}

	We notice that the nonlinear operator $\Ec$ in \eqref{eq:defEc} is in fact a sublinear as the supremum of a family of linear maps on the probability measures.
	In the discrete time gambling house model of Dellacherie \cite{Dellacherie_MaisonJeux}, Dellacherie and Meyer \cite{Dellacherie_Meyer_1988},
    the nonlinear operator can be in fact defined as the supremum of a family of ``nonlinear'' maps on the measures,
    and the time consistency can be deduced by the same arguments.
	We can adapt his discrete time model to our continuous model by introducing a penalty function on the probability measures.
	Suppose that $(\alpha_{s,t})_{s \le t}$ is a family of penalty function where $\alpha_{s,t}: \Pc(\Om) \to \R^+$, we consider the following nonlinear operators:
	\b*
		\Ec_{s, t} [\xi](\om)
		&:=&
		\sup \Big\{ \E^{\P} \big[ -\xi \big] - \alpha_{s,t}(\P) ~ : \P \in \Pc_{s, \om} \Big \}.
	\e*
	
	Assume that for every fixed $t \in \R^+$, $(s, \P) \mapsto \alpha_{s, t}(\P)$ is analytic and the co-cycle condition holds true, i.e. for all $r \le s \le t$,
	\b*
		\alpha_{r,t} (\P)
		&=&
		\alpha_{r,s} (\P) ~+~ \E^{\P} \big[ \alpha_{s,t}(\P_{\om}) \big],
	\e*
	where $(\P_{\om})_{\om \in \Om}$ is a family of r.c.p.d. of $\P$ w.r.t. $\Fc_s$.
	Then under appropriate conditions and by similar arguments as in Theorem \ref{theo:TC}, we can easily obtain the time consistency of $(\Ec_{s,t})_{s \le t}$ of the form
    \b*
        \Ec_{r,t}[\xi] = \Ec_{r,s} \big[ - \Ec_{s,t}[\xi] \big], ~\forall r \le s\le t,~ \xi \in \Ac_{usa}(\Om).
    \e*

	This formulation is closed to the dynamic risk measure proposed by Bion-Nadal \cite{BionNadal2009}, where $\alpha$ is called the minimal penalty function.
	We also notice that in the dynamic risk measure case, $\Ec_{s,t}$ is considered as an operator from $\Ac^U_t$ to $\Ac^U_s$,
	and	$\alpha_{s,t}$ is assumed to satisfy that $\alpha_{s,t}(\P) = \alpha_{s,t}(\Q)$ whenever $\P|_{\Fc_t} = \Q|_{\Fc_t}$.
	An example of the penalty function, which is similar to that considered in Bion-Nadal \cite{BionNadal2012}, is given by
	\b*
		\alpha_{s,t}(\P) &:=& \E^{\P} \Big[ \int_s^t \ell(r, X_r) dr \Big], ~\mbox{for some function}~ \ell : \R^+ \x E \to \R^+.
	\e*

\subsection{Extension on an enlarged space as optimal stopping}
\label{subsec:TC_enlarged_space}

    The above formulation provides a framework to study the optimal control problems.
    Motivated by optimal stopping problems, we introduce an enlarged canonical space $\Omh = \Om \x \R^+$ with current point $\omh = (\om, \theta)$ and historical process $\Yh = (Y_t, \Theta_t)_{t \ge 0}$.
    The process $\Theta$ is defined, for every $\omh = (\om, \theta) \in \Omh$, $\Theta_t(\omh) := \theta 1_{\theta \le t} + \partial 1_{\theta > t}$,
    where $\partial$ is a fictitious time such that $0 + \partial = \partial$.
    By abuse of notation, we also denote $Y(\om)$ in place of $Y(\omh)$, and $\Theta(\theta)$ in place of $\Theta(\omh)$ for $\omh = (\om, \theta) \in \Omh$.

    Let $t, ~\eta, ~\theta \in \R^+$ be such that $\Theta_t(\theta) = \Theta_t(\eta)$, then the concatenation of $\eta$ and $\theta$ at time $t$ is then given by
    $\eta \ox_{t} \theta := \eta 1_{\eta \le t} + \theta 1_{\eta > t}$,
    which is equivalent to $ \Theta_s(\eta\ox_{t} \theta) = \Theta_s(\eta) 1_{s \le t} + \Theta_s(\theta) 1_{s > t}$.
    The enlarged canonical filtration $\widehat \F = (\widehat \Fc_t)_{t \ge 0}$, defined by $\widehat \Fc_t := \sigma(X_s, \Theta_s, s \le t)$, is also generated by the stopping operator $\hat a_t(\omh) := \widehat Y_t(\omh)$.
    For every $(\wh,t) \in \Omh \x \R^+$, denote $\widehat \Dc_{(\wh,t)} := \{ \omh=(\om,\theta) ~: Y_t(\om)=Y_t(\w), \Theta_t( \theta)=\Theta_t( \eta)\}$,
    and $\widehat \Dc_{\wh}^t := \{ \omh=(\om,\theta) ~: X_t(\om)=X_t(\w), \Theta_t( \theta)=\Theta_t( \eta)\}$.
    The concatenation operator for $\wh=(\w, \eta)$ and $\omh = (\om, \theta) \in \widehat \Dc_{(\wh,t)}$ is defined by $\wh\ox_{t} \omh=(\w\ox_{t} \om, \eta\ox_{t} \theta)$.
    Denote by $\widehat{\Oc}$ the optional $\sigma$-field associated with the filtration $\widehat \F$,
    then similarly, $(\wh, t, \omh) \mapsto \wh \ox_t \omh$ restricted on $\big\{ (\wh, t, \omh) ~: \omh \in \widehat \Dc_{\wh}^t \big\}$ is $\widehat{\Oc} \ox \widehat{\Fc}_{\infty}$-measurable.
    We finally notice that $\Theta_{\infty}$, defined by $\Theta_{\infty}(\theta) := \theta$, is in particular a finite $\widehat \F$-stopping time on $\Omh$,
    i.e. $\Theta_{\infty} \in \widehat{\Tc} := \{ \mbox{All finite}~ \widehat{\F} \mbox{-stopping times} \}$.

    Similarly, we can define conditioning and concatenation for probability measures on $\Omh$.
    First, by the same arguments as in item \rmiv of the paragraph ``Canonical filtration'' in Section \ref{subsec:canonical_space},
    $\widehat \Fc_{\hat \tau}$ is countably generated for every finite $\widehat \F$-stopping time $\hat \tau$.
    Then let $\widehat \P$ be a probability measure on $(\Omh, \widehat \Fc_{\infty})$, there is a r.c.p.d. $(\widehat \P_{\omh})_{\omh \in \Omh}$ of $\widehat \P$ w.r.t. $\widehat \Fc_{\hat \tau}$.
    In particular, $\widehat \P_{\omh} (\widehat \Dc_{(\wh, \hat \tau(\wh))}) = 1$ for every $\omh \in \Omh$.
    Next, let $\widehat \P$ be a probability measure on $(\widehat \O, \widehat \Fc_{\hat \tau})$
    and $(\widehat \Q_{\wh})_{\wh \in \Omh}$ be a  kernel probability  from $(\Omh, \widehat \Fc_{\tau})$ to $(\Omh, \widehat \Fc_{\infty})$
    such that $\widehat \Q_{\wh}(\widehat \Dc_{\wh}^{\hat \tau(\wh)}) = 1, ~\forall \wh \in \Omh$,
    then there is a unique concatenated probability measure $\widehat \P \ox_{\hat \tau} \widehat \Q_{\cdot}$ on $(\Omh, \widehat \Fc_{\infty})$, defined by
	\be \label{eq:concatenation_tau_h}
		\widehat \P \ox_{\hat \tau}  \widehat \Q_{\cdot}(A):=\int \widehat \P(d\wh)\int \1_A(\wh \ox_{\hat \tau(\wh)}\omh) \widehat \Q_{\wh}(d\omh).
	\ee

    We then consider a family $(\Pch_{t,\wh})_{(t,\wh) \in \R^+ \x \Omh}$ of collections of probability measures on $\Omh$
    and a family of operators $(\widehat \Ec_{\tauh})_{\tauh \in \widehat{\Tc}}$, on the space $\Ac_{usa}(\Omh)$ of all upper semianalytic functions on $\Omh$, defined by
	\be \label{eq:Ech}
		\widehat \Ec_{\hat \tau} \big[\xi \big](\omh)
		&:=&
		\sup \big\{ \E^{\widehat \P} \big[ \xi \big]\>  ~: \> \widehat \P \in \Pch_{\hat\tau(\omh),\omh} \big\},~~\forall ~\mbox{finite}~\widehat \F-\mbox{stopping time}~ \hat \tau.
	\ee

    \begin{Assumption} \label{assum:stability_h}
        \rmi ({\bf Measurability}) The graph $[[\Pch]]$ is analytic in the Polish space $\R^+ \x \Omh \x \Pc(\Omh)$, where
        \be \label{eq:Pc_graph_h}
            \big[ \big[ \Pch \big] \big]
            ~:=~
            \big \{ (t, \omh, {\widehat \P}) ~: (t,\omh) \in \R^+ \x \Omh,~{\widehat \P} \in \Pch_{t, \omh} \big\} .
        \ee
        \rmii ({\bf Adaptation}) For every $(t,\wh) \in \R^+ \x \Omh$,  $\Pch_{t,\wh} $ is not empty and
        \be \label{eq:Prop_Pc_h}
            \Pch_{t,\wh} = \Pch_{t, [\wh]_t},
            &&
            \widehat \P \big( \widehat \Dc_{(\wh,t)}\big) = 1,~~
            \forall \> \widehat \P \in \Pch_{t,\wh},
        \ee
        \rmiii ({\bf Stability by conditioning and concatenation})
        Let $(t_0,\omh_0) \in \R^+ \x \Omh$, $\hat \tau$ be a stopping time taking value in $[t_0, \infty)$ and $\widehat\P \in \Pch_{t_0,\omh_0}$.
        There is a family of regular conditional probability measures $(\widehat\P_{\omh})_{\omh \in \Omh}$ of $\widehat \P$ w.r.t. $\widehat \Fc_{\hat \tau}$
        such that $\widehat \P_{\omh} \in \Pch_{\hat \tau(\omh), \omh}$ for $\widehat \P$-almost every $\omh \in \Omh$.
        Moreover, let $(\widehat \Q_{\omh})_{\omh \in \Omh}$ be such that $\omh \mapsto \widehat \Q_{\omh}$ is $\widehat \Fc_{\hat \tau}$-measurable
        and $\widehat \Q_{\omh} \in \widehat \Pc_{\hat \tau(\omh), \omh}$ for $\widehat \P$-a.e. $\omh \in \Omh$, then $\widehat \P \ox_{\hat \tau} \widehat \Q_{\cdot} \in \Pch_{t_0, \omh_0}$.
    \end{Assumption}

    By the same arguments as in Theorem \ref{theo:TC}, we have the following time consistence result.
	\begin{Theorem} \label{theo:TCh}
		Let Assumption \ref{assum:stability_h} hold true.
        Then for all finite $\widehat \F$-stopping times $\hat \tau \le \hat \sigma$ and $\xi \in \Ac_{usa}(\Omh)$, $\widehat \Ec_{\hat \tau}[\xi] \in \Ac^U_{\hat \tau}(\Omh) \subset \Ac_{usa}(\Omh)$.
        Further, we have the following time consistence property
		\b*
			  \widehat\Ec_{\hat \tau} [\xi] &=& \widehat \Ec_{\hat \tau} \big[ \widehat \Ec_{\hat \sigma} [\xi] \big] .
		\e*
	\end{Theorem}

\subsection{An abstract dynamic programming principle}
\label{subsec:abstr_dpp}

    The above formulation on the enlarged canonical space $\Omh$ provides a framework to study the optimal control/stopping problem.
    Generally, a stochastic control/stopping term is a controlled stochastic process together with a stopping time.
    The problem consists in maximizing the expected reward value, which depends on the stopped path of the controlled process.
    By considering the distribution on $\Omh$ induced by the controlled process and the stopping time, the problem can be reformulated in the form \eqref{eq:Ech}.
    Moreover, the time consistence property turns to be a dynamic programming principle (DPP) of the control problem.
    Let us provide here an abstract DPP.

    Let $\big(\Pch^0_{t,\xb} \big)_{(t,\xb) \in \R^+ \x \Om }$ be a family of subsets in $\Pc(\Omh)$.
		Similarly in Assumption \ref{assum:stability_h}, we suppose that $\Pch^0_{t,\xb} = \Pch^0_{t,[\xb]_t}$ is nonempty for every $(t,\xb) \in \R^+ \x \Om$ and every $\Ph \in \Pch^0_{t,\xb}$ satisfies that $\Ph \big[ \Theta \ge t,~ X_s = \xb_s,~\forall 0 \le s \le t \big] = 1$.
		In the optimal control/stopping problem context, the above property implies that $\Ph \in \Pch_{t,\xb}$ defines the distribution the controlled process and that of the stopping time after time $t$, while the past path is fixed as $\xb$.
		We suppose further that the graph
		$
				\big[ \big[ \Pch^0 \big]\big]
				:=
				\big\{
					(t,\xb,\Ph) ~: \Ph \in \Pch^0_{t,\xb}
				\big\}
		$
		is analytic in $\R^+ \x \Om \x \Pc(\Omh)$.
		Moreover, the family $\big(\Pch^0_{t,\xb} \big)_{(t,\xb) \in \R^+ \x \Om }$ is stable by conditioning and concatenation in the following sense:
		Let $(t_0, \xb_0) \in \R^+ \x \Om$, $\Ph \in \Pch^0_{t_0, \xb_0}$ and $\tauh$ be a $\widehat{\F}$-stopping time taking value in $[t_0, \infty)$, denote $A_{\tauh} := \{ \omh ~: \Theta_{\infty} > \tauh \}$.
		Then there is a family of r.c.p.d. $(\Ph_{\omh})_{\omh \in \Omh}$ of $\Ph$ w.r.t. $\widehat{\Fc}_{\tauh}$
		such that for $\Ph$-a.e. $\omh = (\om, \theta) \in A_{\tauh}$, $\Ph_{\omh} \in \Pch^0_{\tauh(\omh), \om}$;
		further, let $(\Qh_{\omh})_{\omh \in \Omh}$ be such that $\omh \mapsto \Qh_{\omh}$ is $\widehat{\Fc}_{\tauh}$-measurable,
		$\widehat \Q_{\omh}(\widehat \Dc_{\omh}^{\hat \tau(\omh)}) = 1, ~\forall \omh \in \Omh$
		and $\Qh_{\omh} \in \Pch^0_{\tauh(\omh), \om}$ whenever $\omh \in A_{\tauh}$, then $\Ph \ox_{\tauh} \Qh_{\cdot} \in \Pch^0_{t_0, \xb_0}$.
		
		Let $\Phi : \Omh \to \R^+$ be the positive measurable reward function such that $\Phi(\om,\theta) = \Phi([\om]_{\theta}, \theta)$ for all $(\om, \theta) \in \Omh$,
		the value function of the optimal control/stopping problem is then given by, for all $(t,\xb) \in \R^+ \x \Om$,
		\be \label{eq:defV}
			V(t,\xb)
			&:=&
			\sup_{\Ph \in \Pch^0_{t, \xb}}
			\E^{\Ph} \Big[  \Phi \big( X_{\cdot}, \Theta_{\infty} \big) \Big].
		\ee
		In is clear that $V(t,\xb) = V(t,[\xb]_t)$ since $\Pch^0_{t,\xb} = \Pch^0_{t,[\xb]_t}$, which implies that $V(t,\xb)$ only depends on the past information before $t$ given by $\xb$.

	\begin{Theorem} \label{theo:DPP}
		Let $(\Pch^0_{t,\xb})_{(t,\xb) \in \R^+ \x \Om}$ be the family given above, and $V$ be defined in \eqref{eq:defV}.
		Then $V:\R^+ \x \Om \to \R^+$ is upper semi-analytic and in particular universally measurable.
		Moreover, for every $(t,\xb) \in \R^+ \x \Om$ and every $\widehat{\F}$-stopping time $\tauh$ taking value in $[t,\infty)$, we have the DPP
		\be \label{eq:DPP}
			V(t,\xb)
			&=&
			\sup_{\Ph \in \Pch^0_{t, \xb}} \E^{\Ph}
			\Big[ 1_{\Theta_{\infty} \le \tauh} \Phi \big( X_{\cdot}, \Theta_{\infty} \big)
			~+~
			1_{\Theta_{\infty} > \tauh} V \big( \tauh, [X]_{\tauh} \big)  \Big].
		\ee
	\end{Theorem}
    \proof We use the same arguments as in Theorem \ref{theo:TC}.
    First, the measurability of $V$ is an immediate consequence of the fact that $\big[ \big[ \Pch^0 \big]\big]$ is analytic.
    Then by considering an arbitrary $\Ph \in \Pch^0_{t,\xb}$ as well as its r.c.p.d. $(\Ph_{\omh})_{\omh \in \Omh}$ w.r.t. $\Fch_{\tauh}$, we can easily get that
		\b*
			V(t,\xb)
			& \le &
			\sup_{\Ph \in \Pch^0_{t, \xb}} \E^{\Ph}
			\Big[ 1_{\Theta_{\infty} \le \tauh} \Phi \big( X_{\cdot}, \Theta_{\infty} \big)
			~+~
			1_{\Theta_{\infty} > \tauh} V \big( \tauh, [X]_{\tauh} \big)  \Big].
		\e*
		For the reverse inequality, we let $\eps > 0$, and choose a family $(\Qh^{\eps}_{\wh})_{\wh \in A_{\tauh}}$ such that $\wh \mapsto \Qh^{\eps}_{\wh}$ restricted on $A_{\tauh}$ is $\Fch_{\tauh}$-measurable,
        and $\E^{\Qh^{\eps}_{\wh}} \big[ \Phi \big( X_{\cdot}, \Theta_{\infty} \big) \big] \ge V_{\eps}(\tau(\wh), \wh)$ for every $\wh \in A_{\tauh}$,
        where $V_{\eps}(t,\xb) := (V(t, \xb) - \eps)1_{V(t,\xb)<\infty} + \frac{1}{\eps} 1_{V(t,\xb) = \infty}$.
		Next, we complete the family $(\Qh^{\eps}_{\wh})_{\wh \in A_{\tau}}$ in an arbitrary but measurable way such that $\Ph \ox_{\tauh} \Qh^{\eps}_{\cdot} \in \Pch^0_{t,\xb}$. It follows that
		\b*
			V(t,\xb)
			& \ge &
			\sup_{\Ph \in \Pch^0_{t, \xb}} \E^{\Ph}
			\Big[ 1_{\Theta_{\infty} \le \tauh} \Phi \big( X, \Theta_{\infty} \big)
			~+~
			1_{\Theta_{\infty} > \tauh} V_{\eps} \big( \tauh, [X]_{\tauh} \big)  \Big],
		\e*
		which completes the proof by the arbitrariness of $\eps > 0$.
    \qed

\section{Conclusion}

	We gave a brief introduction to the capacity theory of Choquet.
	Following Dellacherie \cite{Dellacherie_1972}, 
	we showed how to derive a measurable selection theorem using the projection capacity on the product space of an abstract measurable space and a topological space.
	It is classical to use measurable selection techniques to deduce a dynamic programming for discrete time optimization problems. 
	We then also proposed an abstract framework for the dynamic programming principle, or equivalently the time consistency property for a class of continuous time optimization problems.
	In our acompanying paper \cite{ElKaroui_Tan2}, we shall	show that this framework is convient to study the general stochastic control/stopping problem.


\end{document}